\documentclass{amsart}

\usepackage{amsmath}
\usepackage{amsfonts}
\usepackage{amssymb,enumerate}
\usepackage{amsthm}
\usepackage[all]{xy}

\newtheorem{lem}{Lemma}[section]
\newtheorem{cor}[lem]{Corollary}
\newtheorem{prop}[lem]{Proposition}
\newtheorem{thm}[lem]{Theorem}

\newtheorem{Defn}[lem]{Definition}
\newtheorem{Ex}[lem]{Example}
\newtheorem{Question}[lem]{Question}
\newtheorem{Property}[lem]{Property}
\newtheorem{Properties}[lem]{Properties}
\newtheorem{Subprops}{}[lem]
\newtheorem{Para}[lem]{}

\newenvironment{defn}{\begin{Defn}\rm}{\end{Defn}}
\newenvironment{ex}{\begin{Ex}\rm}{\end{Ex}}

\newenvironment{property}{\begin{Property}\rm}{\end{Property}}
\newenvironment{properties}{\begin{Properties}\rm}{\end{Properties}}
\newenvironment{subprops}{\begin{Subprops}\rm}{\end{Subprops}}
\newenvironment{para}{\begin{Para}\rm}{\end{Para}}

\theoremstyle{definition}

\newcommand{\wt}{\widetilde}
\newcommand{\te}{\theta}
\newcommand{\comp}[1]{\widehat{#1}}
\newcommand{\ideal}[1]{\mathfrak{#1}}
\newcommand{\m}{\ideal{m}}
\newcommand{\n}{\ideal{n}}
\newcommand{\p}{\ideal{p}}
\newcommand{\q}{\ideal{q}}

\newcommand{\gdim}{\mathrm{G\text{-}dim}}
\newcommand{\gfd}{\mathrm{Gfd}}   
\newcommand{\pdim}{\mathrm{pd}} 
\newcommand{\rfd}{\operatorname{Rfd}}
\newcommand{\shift}{{\scriptstyle{\Sigma}}}
       
\newcommand{\ext}{\mathrm{Ext}}

\newcommand{\depth}{\mathrm{depth}}     
\newcommand{\rank}{\mathrm{rank}}       

\newcommand{\injdim}{\mathrm{id}}       
\newcommand{\rhom}{\mathbf{R}\mathrm{Hom}}      
\newcommand{\lotimes}{\otimes^{\mathbf{L}}}
\newcommand{\lch}{\mathbf{R}\Gamma}
\newcommand{\amp}{\mathrm{amp}}
\newcommand{\HH}{\mathrm{H}}
\newcommand{\Hom}{\mathrm{Hom}} 
\newcommand{\embdim}{\mathrm{edim}}
\newcommand{\coker}{\mathrm{Coker}}
\newcommand{\fd}{\mathrm{fd}}
\newcommand{\fdim}{\mathrm{fd}}
\newcommand{\spec}{\mathrm{Spec}}
\newcommand{\soc}{\mathrm{Soc}}
\newcommand{\crs}{\mathrm{crs}}
\newcommand{\s}{\mathrm{s}}
\newcommand{\x}{\mathbf{x}}
\newcommand{\tor}{\mathrm{Tor}}
\newcommand{\vf}{\varphi}
\newcommand{\A}{\mathcal{A}}

\begin{document}

\bibliographystyle{amsplain}

\author{Srikanth Iyengar} \address{Mathematics Department, University of Missouri,
  Columbia, MO 65211 USA} 
\curraddr{918 Oldfather Hall, Department of Mathematics, University of Nebraska,
Lincoln, NE 68588}
\email{siyengar@math.unl.edu}

\author{Sean Sather-Wagstaff} \address{Department of Mathematics, University of Illinois,
  273 Altgeld Hall, 1409 West Green Street, Urbana, IL, 61801 USA}
\email{ssather@math.uiuc.edu}

\thanks{This research was conducted while S.I. was partly supported by a grant from the
  NSF, and S.S.-W. was an NSF Mathematical Sciences Postdoctoral Research Fellow. Part of
  the work was done at MSRI during the Spring semester of 2003, when the authors were
  participating in the program in Commutative Algebra.}

\title[G-dimension over local homomorphisms]{G-dimension over local homomorphisms.
  Applications to the {F}robenius endomorphism}



\begin{abstract}
  We develop a theory of G-dimension over local homomorphisms which encompasses the
  classical theory of G-dimension for finitely generated modules over local rings.  As an
  application, we prove that a local ring $R$ of characteristic $p$ is Gorenstein if and
  only if it possesses a nonzero finitely generated module of finite projective dimension
  that has finite G-dimension when considered as an $R$-module via some power of the
  Frobenius endomorphism of $R$.  We also prove results that track the behavior of
  Gorenstein properties of local homomorphisms under composition and decomposition.
\end{abstract}
\maketitle

\section{Introduction}

The main goal of this article is to develop a theory of Gorenstein dimension over local
homomorphisms.  More precisely, given a local homomorphism $\vf\colon R\to S$, to each
finitely generated (in short: finite) $S$-module $M$, we attach an invariant
$\gdim_{\vf}(M)$, called the \emph{G-dimension of $M$ over $\vf$}.  This invariant is
defined using the technology of Cohen factorizations, developed by Avramov, Foxby, and
B.~Herzog~\cite{avramov:solh}.  The reader can refer to Section~\ref{sec:gd} for the
details.  When $M$ happens to be finite over $R$, for instance when $\vf=\mathrm{id}_R$,
this coincides with the G-dimension of $M$ over $R$ as defined by Auslander and
Bridger~\cite{auslander:smt}; this is contained in Corollary~\ref{cor:vf-is-finite}.

One of the guiding examples for this work is the Frobenius map $\vf\colon R\to R$, given
by $x\mapsto x^p$, where $R$ is a local ring of positive prime characteristic $p$.  Since
$\vf$ is a ring homomorphism, so is $\vf^n$ for each integer $n>0$, and hence one can view
$R$ as a left module over itself via $\vf^n$.  Denote this $R$-module ${}^{\vf^n}\!\! R$.
Like in the case of the residue field, it is known that certain homological properties of
${}^{\vf^n}\!\!R$ determine and are determined by ring-theoretic properties of $R$.
Consider, for instance, regularity. The Auslander-Buchsbaum-Serre theorem says that a
local ring is regular if and only if its residue field has finite projective dimension.
Compare this with the fact that, when $R$ has characteristic $p$, it is regular if and
only if the flat dimension of ${}^{\vf^n}\!\! R$ is finite for some $n\geq 1$; this is
proved by Kunz~\cite[(2.1)]{kunz:corlrocp} and Rodicio~\cite[(2)]{rodicio:oaroa}.  This
result may be reformulated as: the local ring $R$ is regular if and only if $\pdim(\vf^n)$
is finite for some integer $n\geq 1$.  Here, given any local homomorphism $\vf\colon R\to
S$, we write $\pdim_{\vf}(\text{-})$ for the \emph{projective dimension over $\vf$}, which
is also defined via Cohen factorizations, and $\pdim(\vf)=\pdim_{\vf}(S)$; see
Section~\ref{sec:pd}.

A key contribution of this paper, Theorem A below, is a similar characterization of the
Gorenstein property for $R$. It is contained in Theorem~\ref{thm:gdim(f)-finite} and is
analogous to a classical result of Auslander and Bridger: for any local ring, the residue
field has finite G-dimension if and only if the ring is Gorenstein.

\medskip

\noindent \textbf{Theorem A.}  
\textit{Let $R$ be a local ring of positive prime characteristic $p$ and $\vf$ its
  Frobenius endomorphism.  The following conditions are equivalent.
\begin{enumerate}[{\quad\rm(a)}]
\item[(a)] The ring $R$ is Gorenstein.
\item[(b)] $\gdim(\vf^n)$ is finite for some integer $n\geq 1$.
\item[(c)] There exists a nonzero finite $R$-module $P$ of finite projective dimension and
  an integer $n\geq 1$ such that $\gdim_{\vf^n}(P)$ is finite.
\end{enumerate}}

In the statement, $\gdim(\vf^n)=\gdim_{\vf^n}(R)$.  In the special case where $\vf$ is
module-finite, the equivalence of conditions (a) and (b) coincides with a recent result of
Takahashi and Yoshino~\cite[(3.1)]{takahashi:ccmlrbfm}.  These are related also to a
theorem of Goto~\cite[(1.1)]{goto:aponlrocp}.

The bulk of the article is dedicated to a systematic investigation of the invariant
$\gdim_{\vf}(\text{-})$.  Some of the results obtained extend those concerning the
classical invariant $\gdim_R(\text{-})$.  Others are new even when specialized to the
absolute situation.  The ensuing theorem is one such.  It is comparable
to~\cite[(3.2)]{foxby:daafuc}, which can be souped up to: if $\pdim_{\sigma}(P)$ is
finite, then $\pdim_{\sigma\vf}(P)=\pdim(\vf)+\pdim_{\sigma}(P)$; see
Theorem~\ref{thm:stab-pd-tensor} for a further enhancement.  \medskip

\noindent \textbf{Theorem B.}  
\textit{Let $\vf\colon R\to S$ and $\sigma\colon S\to T$ be local homomorphisms, and let
  $P$ be a nonzero finite $T$-module.  If $\pdim_{\sigma}(P)$ is finite, then
\[ \gdim_{\sigma\vf}(P)=\gdim(\vf)+\pdim_{\sigma}(P).\]
In particular, $\gdim_{\sigma\vf}(P)$ and $\gdim(\vf)$ are simultaneously finite.}

\medskip

This result is subsumed by Theorem~\ref{thm:stability-tensor}.  The special case $P=T$,
spelled out in Theorem~\ref{thm:gdim-compose}, may be viewed as a
composition-decomposition theorem for maps of finite G-dimension.  It is expected that the
composition part of the result holds even when $\gdim(\sigma)$ is
finite~\cite[(4.8)]{avramov:rhafgd}.  However, as Example~\ref{ex:not-weaken}
demonstrates, the decomposition part cannot extend to that generality.

Theorem B and its counterpart for projective dimension are crucial ingredients in the
following theorem that generalizes~\cite[(4.6.c)]{avramov:lgh}
and~\cite[(8.8)]{avramov:rhafgd} proved by Avramov and Foxby.

\medskip

\noindent \textbf{Theorem C.}  
\textit{Let $\vf\colon (R,\m)\to (S,\n)$ and $\sigma\colon (S,\n)\to (T,\p)$ be local
  homomorphisms with $\pdim(\sigma)$ finite.  If $\sigma\vf$ is (quasi-)Gorenstein at
  $\p$, then $\vf$ is (quasi-)Gorenstein at $\n$ and $\sigma$ is Gorenstein at $\p$.}

\medskip

This result coincides with Theorem~\ref{thm:fin-fd-descent}.  Section~\ref{sec:descent}
contains other results of this flavor.  It is worth remarking that there is an analogue of
Theorem C for complete intersection homomorphisms, due to
Avramov~\cite[(5.7)]{avramov:lci}.

It turns out that the \emph{finiteness} of $\gdim_{\vf}(M)$ depends only on the $R$-module
structure on $M$, although its value depends on $\vf$; this is the content of
Theorem~\ref{thm:indep-of-vf} and Example~\ref{ex:not-equal}. One way to understand this
result would be to compare Gorenstein dimension over $\vf$ to various extensions of the
classical G-dimension to $R$-modules that may not be finite.  The last section deals with
this problem, where Theorem~\ref{thm:gdvsgfd} contains the following result; in it
$\gfd_R(M)$ is the Gorenstein flat dimension of $M$ over $R$.

\medskip

\noindent\textbf{Theorem D.}
\textit{ Assume $R$ is a quotient of a Gorenstein ring and let $\vf\colon R\to S$ be a
  local homomorphism. For each finite $S$-module $M$, one has
\[ 
\gfd_R(M) - \embdim(\vf) \leq \gdim_{\vf}(M) \leq \gfd_R(M)\,.
 \]
 In particular, $\gdim_{\vf}(M)$ is finite if and only if $\gfd_R(M)$ is finite.}

\medskip

Foxby, in an unpublished manuscript, has obtained the same conclusion assuming only that
the formal fibres of $R$ are Gorenstein. Specializing $X$ to $S$ yields that $\gdim(\vf)$
and $\gfd_RS$ are simultaneously finite.  This last result was proved also by Christensen,
Frankild, and Holm \cite[(5.2)]{christensen:new}, and our proof of Theorem D draws heavily
on their work.

En route to the proof of Theorem D, we obtain results on G-flat dimension that are of
independent interest; notably, the following Auslander-Buchsbaum type formula for the
depth of a module of finite G-flat dimension.  It is contained in
Theorem~\ref{prop:gfd-depth}.

\medskip

\noindent \textbf{Theorem E.}  
\textit{Let $(R,\m,k)$ be a local ring and $E$ the injective hull of $k$. If $M$ is an
  $R$-module with $\gfd_R(M)$ finite, then}
\[
\depth_R(M) = \depth R - \sup(E\lotimes_RM)\,.
\]

In the preceding discussion, we have focused on modules.  However, most of our results are
stated and proved for complexes of $R$-modules.  This is often convenient and sometimes
necessary, as is the case in Theorem~\ref{thm:stability-tensor}. Section \ref{sec:back} is
mainly a catalogue of standard notions and techniques from the homological algebra of
complexes required in this work; most of them can be found in Foxby's
notes~\cite{foxby:hacr} or Christensen's monograph~\cite{christensen:gd}.

\section{Background} \label{sec:back}

Let $R$ be a commutative Noetherian ring.  A complex of $R$-modules is a sequence of
$R$-module homomorphisms
\[ 
X= \cdots\xrightarrow{\partial_{i+1}}X_i\xrightarrow{\partial_{i}}X_{i-1}
\xrightarrow{\partial_{i-1}}\cdots \] such that $\partial_i\partial_{i+1}=0$ for all $i$.
The \textit{supremum}, the \textit{infimum}, and the \textit{amplitude} of a complex $X$
are defined by the following formulas:
\begin{align*}
  \sup(X) & = \sup\{i\mid \HH_i(X)\neq 0\} \\
  \inf(X) & = \inf\{i\mid \HH_i(X)\neq 0\} \\
  \amp(X) & = \sup(X)-\inf(X).
\end{align*}
Note that $\amp(X)=-\infty$ if and only if $\HH(X)=0$.  The complex $X$ is
\textit{homologically bounded} if $\amp(X)<\infty$, and it is \textit{homologically
  degreewise finite} if $\HH(X)$ is degreewise finite.  When $\HH(X)$ is both degreewise
finite and bounded we say that $X$ is \emph{homologically finite}.

Let $X$ and $Y$ be complexes of $R$-modules.  As is standard, we write $X\lotimes_R Y$ for
the derived tensor product of $X$ and $Y$, and $\rhom_R(X,Y)$ for the derived
homomorphisms from $X$ to $Y$.  The symbol ``$\simeq$'' denotes an isomorphism in the
derived category.  For details on derived categories and derived functors, the reader may
refer to the classics, Hartshorne~\cite{hartshorne:rad} and
Verdier~\cite{verdier}, or, for a more recent treatment, to Gelfand and
Manin~\cite{gelfand:moha}.

Let $X$ be a homologically bounded complex of $R$-modules.  A \textit{projective}
\textit{resolution} of $X$ is a complex of projective modules $P$ with $P_i=0$ for $i\ll
0$ and equipped with an isomorphism $P\simeq X$.  Such resolutions exist and can be chosen
to be degreewise finite when $X$ is homologically finite.  The \textit{projective}
\textit{dimension} of $X$ is
\[ 
\pdim_R(X):=\inf\{\sup\{n\mid P_n\neq 0\}\mid \text{$P$ a projective resolution of $X$}\}.
\] 
Thus, if $\HH(X)=0$, then $\pdim_R(X)$ is $-\infty$, and hence it is not finite.  Flat
resolutions and injective resolutions, and the corresponding dimensions $\fdim_R(X)$ and
$\injdim_R(X)$, are defined analogously.

The focus of this paper is G-dimension for complexes.  In the next few paragraphs, we
recall its definition and certain crucial results that allow one to come to grips with it.

\begin{para}
  A finite $R$-module $G$ is \textit{totally reflexive} if
\begin{enumerate}[{\quad\rm(a)}]
\item $\ext^i_R(G,R)=0$ for all $i>0$;
\item $\ext^i_R(G^*,R)=0$ for all $i>0$, where $(\text{-})^*$ denotes
  $\Hom_R(\text{-},R)$; and
\item the canonical map $G\to G^{**}$ is bijective.
\end{enumerate}

Let $X$ be a homologically finite complex of $R$-modules.  A \textit{G-resolution} of $X$
is an isomorphism $G\simeq X$ where $G$ is complex of totally reflexive modules with
$G_i=0$ for $i\ll 0$.  A degreewise finite projective resolution of $X$ is also a
G-resolution, since every finite projective module is totally reflexive.  The
\textit{G-dimension} of $X$ is
\[ \gdim_R(X):=\inf\{\sup\{n\mid G_n\neq 0\}\mid \text{$G$ is a  
  G-resolution of $X$}\}.\]
\end{para}

The following paragraphs describe alternative, and often more convenient, ways to detect
when a complex has finite G-dimension.

\begin{para} \label{para:reflexive}
  A homologically finite complex $X$ of $R$-modules is \textit{reflexive} if
\begin{enumerate}[{\quad\rm(a)}]
\item $\rhom_R(X,R)$ is homologically bounded; and
\item the canonical biduality morphism below is an isomorphism
\[ \delta^R_X\colon X\to\rhom_R(\rhom_R(X,R),R). \]
\end{enumerate}
\end{para}

This notion is relevant to this article because of the next result, based on an
unpublished work of Foxby; see~\cite[(2.3.8)]{christensen:gd}
and~\cite[(2.7)]{yassemi:gd}.

\begin{para} \label{para:reflexive-2}
  \emph{The complex $X$ is reflexive if and only if\, $\gdim_R(X)<\infty$.  When $X$ is
    reflexive, $\gdim_R(X)=-\inf(\rhom_R(X,R))$.}
\end{para}

Using this characterization, it is easy to verify the base change formula below;
Christensen~\cite[(5.11)]{christensen:sdc} has established a much stronger statement.

\begin{para} 
\label{lem:fflat}
\emph{Let $R\to S$ be a flat local homomorphism and $X$ a homologically finite complex of
  $R$-modules.  Then $\gdim_R(X)=\gdim_S(X\otimes_R S)$.}
\end{para}

Henceforth, $R$ is a local ring, where ``local'' means ``local and Noetherian''.

\begin{para} \label{para:auslander}
  A \textit{dualizing complex} for $R$ is a homologically finite complex of $R$-modules
  $D$ of finite injective dimension such that the natural map $R\to\rhom_R(D,D)$ is an
  isomorphism.  When $R$ is a homomorphic image of a Gorenstein ring, for example, when
  $R$ is complete, it possesses a dualizing complex.
  
  Assume that $R$ possesses a dualizing complex $D$.  The \textit{Auslander category of
    $R$}, denoted $\A(R)$, is the full subcategory of the derived category of $R$ whose
  objects are the homologically bounded complexes $X$ such that
\begin{enumerate}[{\quad\rm(a)}]
\item $D\lotimes_R X$ is homologically bounded; and
\item the canonical morphism below is an isomorphism
\[ \gamma_X\colon X\to\rhom_R(D,D\lotimes_R X). \] 
\end{enumerate}
\end{para}

It should be emphasized that a complex can be in the Auslander category of $R$ without
being homologically finite.  Those that are homologically finite are identified by the
following result; see~\cite[(3.1.10)]{christensen:gd} for a proof.

\begin{para} \label{para:auslander-2}
  \emph{Let $X$ be a homologically finite complex. Then $X$ is in $\A(R)$ if and only if
    $\gdim_R(X)<\infty$.}
\end{para}

The various homological dimensions are related to another invariant: depth.

\begin{para} \label{para:depth}
  Let $K$ be the Koszul complex on a generating sequence of length $n$ for the maximal
  ideal of $R$.  The \textit{depth} of $X$ is defined to be
\[ \depth_R(X)=n-\sup (K\otimes_R X).\]
It is independent of the choice of generating sequence and may be calculated via the
vanishing of appropriate local cohomology or $\ext$-modules~\cite[(2.1)]{foxby:daafuc}.
\end{para}

For the basic properties of depth, we refer to~\cite{foxby:daafuc}.  However, there seems
to be no available reference for the following result.

\begin{lem} \label{lem:depth-eq}
  Let $\vf\colon R\to S$ be a local homomorphism and $X$ a complex of $S$-modules. If \,
  $\HH(X)$ is degreewise finite over $R$, then $\depth_S(X)=\depth_R(X)$.
\end{lem}

\begin{proof}
  Let $K$ denote the Koszul complex on a set of $n$ generators for the maximal ideal of
  $R$.  Note that $\pdim_S(K\otimes_R S)=n$.  Thus
\begin{align*} 
  \depth_S(K\otimes_R X)&=\depth_S((K\otimes_R S)\otimes_S X)\\
  &=\depth_S(X)-\pdim_S(K\otimes_R S)\\
  &=\depth_S(X)-n
\end{align*}
where the second equality is by the Auslander-Buchsbaum
formula~\cite[(2.4)]{foxby:daafuc}.  Now, $\HH(K\otimes_R X)$ is degreewise finite over
$R$ and is annihilated by the maximal ideal of $R$; see, for instance,
\cite[(1.2)]{iyengar:dfcait}.  Hence, each $\HH_i(K\otimes_R X)$ has finite length over
$R$, and, therefore, over $S$.  In particular, by~\cite[(2.7)]{foxby:daafuc} one has
$\depth_S(K\otimes_R X)=-\sup(K\otimes_R X)$.  Combining this with the displayed formulas
above yields that
\[ \depth_S(X)=n+\depth_S(K\otimes_R X)=n-\sup(K\otimes_R X)=\depth_R(X). \]
This is the desired equality.
\end{proof}

In \cite[(3.1)]{foxby:daafuc}, Foxby and Iyengar extend Iversen's Amplitude Inequality; we
require a slight reformulation of their result.

\begin{thm} \label{thm:GAI}
  Let $S$ be a local ring and let $P$ be a homologically finite complex of $S$-modules
  with $\pdim_S(P)$ finite.  For each homologically degreewise finite complex $X$ of
  $S$-modules, one has
\[ \amp (X) \leq\amp(X\lotimes_S P) \leq\amp(X)+\pdim_S(P)-\inf(P).\]
In particular, $\amp (X)$ is finite if and only if $\amp(X\lotimes_S P)$ is finite.
\end{thm}

\begin{proof}
  The inequality on the left is contained in~\cite[(3.1)]{foxby:daafuc}, while the one on
  the right is by~\cite[(7.28), (8.17)]{foxby:hacr}.
\end{proof}

Here is a corollary; one can give a direct proof when the map $\alpha$ is between
complexes that are homologically bounded to the right.

\begin{prop} \label{prop:tool1}
  Let $S$ be a local ring, $P$ a homologically finite complex of $S$-modules with
  $\pdim_S(P)$ finite, and let $\alpha$ be a morphism of homologically degreewise finite
  complexes.  Then $\alpha$ is an isomorphism if and only if the induced map
  $\alpha\lotimes_S P$ is an isomorphism.
\end{prop}

\begin{proof}
  Let $\mathrm{C}(\alpha)$ and $\mathrm{C}(\alpha\lotimes_S P)$ denote the mapping cones
  of $\alpha$ and $\alpha\lotimes_S P$, respectively.  The homology long exact sequence
  arising from mapping cones yields that $\HH(\mathrm{C}(\alpha))$ is degreewise finite.
  Observe that $\mathrm{C}(\alpha\lotimes_S P)=\mathrm{C}(\alpha)\lotimes_S P$.  By the
  previous theorem, $\HH(\mathrm{C}(\alpha))=0$ if and only if
  $\HH(\mathrm{C}(\alpha\lotimes_S P))=0$.
\end{proof}

It is well known that the derived tensor product of two homologically finite complexes is
homologically finite when one of them has finite projective dimension.  In the sequel we
require the following slightly more general result, contained in
\cite[(4.7.F)]{avramov:hdouc}. The proof is short and simple, and bears repetition.

\begin{lem} \label{lem:back-finite}
  Let $\sigma\colon S\to T$ be a local homomorphism and let $X$ and $P$ be homologically
  finite complexes of modules over $S$ and $T$, respectively.  If $\fd_S(P)$ is finite,
  then the complex of $T$-modules $X\lotimes_S P$ is homologically finite.
\end{lem}

\begin{proof}
  Replacing $X$ by a soft truncation, one may assume that $X$ is bounded; see, for
  example, \cite[p.~165]{christensen:gd}.  With $F$ a bounded flat resolution of $P$ over
  $S$, the complex $X\lotimes_S P$ is isomorphic to $X\otimes_S F$, which is bounded.
  Thus, $X\lotimes_S P$ is homologically bounded.  As to its degreewise finiteness: let
  $Y$ and $Q$ be minimal free resolutions of $X$ and $P$ over $S$ and $T$, respectively.
  Then $X\lotimes_S P$ is isomorphic to $Y\otimes_S Q$, which is a complex of finite
  $T$-modules.  Therefore, the same is true of its homology, since $T$ is Noetherian.
\end{proof}

\section{G-dimension over a local homomorphism} \label{sec:gd}

In this section we introduce the G-dimension over a local homomorphism and document some
of its basic properties.  We begin by recalling the construction of Cohen factorizations
of local homomorphisms as introduced by Avramov, Foxby, and B.~Herzog~\cite{avramov:solh}.

\begin{para}  \label{def:factor} 
  Given a local homomorphism $\vf\colon (R,\m)\to (S,\n)$, the \textit{embedding
    dimension} and \emph{depth of $\vf$} are
\[ \embdim(\vf):=\embdim(S/\m S)\qquad\text{and}\qquad
\depth(\vf):=\depth(S)-\depth(R).\] A \textit{regular} (respectively, \textit{Gorenstein})
\textit{factorization} of $\vf$ is a diagram of local homomorphisms,
$R\xrightarrow{\Dot{\vf}}R'\xrightarrow{\vf'}S$, where $\vf=\vf'\Dot{\vf}$, with
$\Dot{\vf}$ flat, the closed fibre $R'/\m R'$ regular (respectively, Gorenstein) and
$\vf'\colon R'\to S$ surjective.

Let $\comp{S}$ denote the completion of $S$ at its maximal ideal and $\iota\colon
S\to\comp{S}$ be the canonical inclusion.  By~\cite[(1.1)]{avramov:solh} the composition
$\grave{\vf}=\iota\vf$ admits a regular factorization $R\to R'\to \comp{S}$ with $R'$
complete.  Such a regular factorization is said to be a \textit{Cohen factorization} of
$\grave{\vf}$.
\end{para}

The result below is analogous to~\cite[(4.3)]{avramov:rhafgd}.  Here, and elsewhere, we
write $\comp{X}$ for $X\otimes_S \comp{S}$ when $X$ is a complex of $S$-modules.

\begin{thm}  \label{thm:likeAF4.3}
  Let $\vf\colon R\to S$ be a local homomorphism and $X$ a homologically finite complex of
  $S$-modules.  If $R\stackrel{\Dot{\vf_1}}{\to} R_1\stackrel{\vf_1'}{\to} \comp{S}$ and
  $R\stackrel{\Dot{\vf_2}}{\to} R_2\stackrel{\vf_2'}{\to} \comp{S}$ are Cohen
  factorizations of $\grave{\vf}$, then
\[ \gdim_{R_1}(\comp{X})-\embdim(\Dot{\vf_1})
= \gdim_{R_2}(\comp{X})-\embdim(\Dot{\vf_2}).\]
\end{thm}
   
\begin{proof}
  Theorem~\cite[(1.2)]{avramov:solh} provides a commutative diagram
\[ \xymatrix{
  & R_1 \ar[dr]^{\vf_1'} \\
  R \ar[ur]^{\Dot{\vf}_1} \ar[r]^{\Dot{\vf}} \ar[dr]_{\Dot{\vf}_2} &
  R' \ar[u]^{v_1} \ar[r]^{\vf'} \ar[d]_{v_2} & \comp{S} \\
  & R_2\ar[ur]_{\vf_2'} } \] where $\vf' \Dot{\vf}$ is a third Cohen factorization of
$\grave{\vf}$, and each $v_i$ is surjective with kernel generated by an $R'$-regular
sequence whose elements are linearly independent over $R'/\m'$ in $\m'/((\m')^2+\m R')$.
Here $\m$ and $\m'$ are the maximal ideals of $R$ and $R'$, respectively.  Let $c_i$
denote the length of a regular sequence generating $\ker(v_i)$.  For $i=1,2$ one has that
\begin{align*}
  \gdim_{R'}(\comp{X})-\embdim(\Dot{\vf})
  & = [\gdim_{R_i}(\comp{X})+c_i]-[\embdim(R_i/\m R_i)+c_i] \\
  & = \gdim_{R_i}(\comp{X})-\embdim(\Dot{\vf_i})
\end{align*}
where~\cite[(2.3.12)]{christensen:gd} gives the first equality.  This gives the desired
result.
\end{proof}

\begin{defn}  \label{def:gdim}
  Let $\vf\colon R\to S$ be a local homomorphism and $X$ a homologically finite complex of
  $S$-modules.  Let $R\xrightarrow{\Dot{\vf}} R'\xrightarrow{\vf'} \comp{S}$ be a Cohen
  factorization of $\grave{\vf}$.  The \textit{G-dimension of $X$ over $\vf$} is the
  quantity
\[ \gdim_{\vf}(X):= 
\gdim_{R'}(\comp{X})-\embdim(\Dot{\vf}). \] Theorem~\ref{thm:likeAF4.3} shows that
$\gdim_{\vf}(X)$ does not depend on the choice of Cohen factorization.  Note that
$\gdim_{\vf}(X)\in\{-\infty\}\cup\mathbb{Z}\cup\{\infty\}$, and also that
$\gdim_{\vf}(X)=-\infty$ if and only if $\HH(X)=0$.

The \textit{G-dimension of $\vf$} is defined to be
\[ \gdim(\vf):=\gdim_{\vf}(S).\]
It is clear from the definitions that the corresponding notion of the finiteness of
$\gdim(\vf)$ agrees with that in~\cite{avramov:rhafgd}.
\end{defn}

Here are some properties of the $\gdim_{\vf}(\text{-})$.

\begin{properties} \label{props}
  Fix a local homomorphism $\vf\colon R\to S$, a Cohen factorization $R\to R'\to \comp{S}$
  of $\grave{\vf}$, and a homologically finite complex $X$ of $S$-modules.

\begin{subprops} \label{subprop:1}
  Let $\comp{\vf}\colon\comp{R}\to\comp{S}$ denote the map induced on completions.  One
  has
\[ \gdim_{\vf}(X)=\gdim_{\grave{\vf}}(\comp{X})
=\gdim_{\comp{\vf}}(\comp{X}).
 \] More generally, let $I$ and $J$ be proper ideals of $R$
and $S$, respectively, with $IS\subseteq J$, and let $\wt{R}$ and $\wt{S}$ denote the
respective completions.  With $\wt{\vf}\colon\wt{R}\to\wt{S}$ the induced map, one has
\[ 
\gdim_{\vf}(X)=\gdim_{\wt{\vf}}(\wt{S}\otimes_S X). 
\]
This is because the completion of $\wt{\vf}$ at the maximal ideal of $\wt{S}$ is
$\comp{\vf}$.
\end{subprops}

\begin{subprops} \label{subprop:3}
  If $X\simeq X'\oplus X''$, then
  $\gdim_{\vf}(X)=\max\{\gdim_{\vf}(X'),\gdim_{\vf}(X'')\}$; this follows from the
  corresponding property of the classical G-dimension.
\end{subprops}

\begin{subprops} \label{subprop:5}
  If $\vf$ has a regular factorization $R\xrightarrow{\vf_1} R_1\xrightarrow{\vf'} S$,
  then
\[ 
\gdim_{\vf}(X)=\gdim_{R_1}(X)-\embdim(\vf_1), 
\]
because the diagram $R\xrightarrow{\grave{\vf_1}} \comp{R_1}\xrightarrow{\comp{\vf'}}
\comp{S}$ is a Cohen factorization of $\grave{\vf}$.
\end{subprops}

\begin{subprops} \label{subprop:4}
  If $\vf$ is surjective, then $R\xrightarrow{=}R\xrightarrow{\vf} S$ is a regular
  factorization, so 
\[ \gdim_{\vf}(X)=\gdim_R(X). \]
Corollary~\ref{cor:vf-is-finite} below generalizes this to the case when $\HH(X)$ is
finite over $R$.
\end{subprops}
\end{properties}

The following theorem is an extension of the Auslander-Bridger formula, which is the
special case $\vf=\mathrm{id}_R$.

\begin{thm} \label{thm:AB}
  Let $\vf\colon R\to S$ be a local homomorphism and $X$ a homologically finite complex of
  $S$-modules.  If $\gdim_{\vf}(X)<\infty$, then
\[ 
\gdim_{\vf}(X)=\depth(R)-\depth_S(X)
 \] 
\end{thm}

\begin{proof}
  Let $R\to R'\xrightarrow{\vf'}\comp{S}$ be a Cohen factorization of $\grave{\vf}$, and
  let $\m$ be the maximal ideal of $R$.  The classical Auslander-Bridger formula gives the
  first of the following equalities; the flatness of $R\to R'$ and the surjectivity of
  $\vf'$ imply the second; the regularity of $R'/\m R'$ yields the third.
\begin{align*}
  \gdim_{R'}(\comp{X})
  & = \depth(R') - \depth_{R'}(\comp{X}) \\
  & = [\depth(R)+\depth(R'/\m R')] - \depth_{\comp{S}}(\comp{X}) \\
  & = \depth(R)-\depth_S(X)+\embdim(R'/\m R')
\end{align*}
This gives the desired equality.
\end{proof}

As in the classical case, described in~\ref{para:auslander-2}, when $R$ has a dualizing
complex one can detect finiteness of $\gdim_{\vf}(\text{-})$ in terms of membership in the
Auslander category of $R$.

\begin{prop} \label{prop:likeAF4.3}
  Let $\vf\colon R\to S$ be a local homomorphism and $R\xrightarrow{\dot{\vf}}
  R'\xrightarrow{\vf'} \comp{S}$ a Cohen factorization of $\grave{\vf}$.
  The following conditions are equivalent for each homologically finite complex $X$ of
  $S$-modules.
\begin{enumerate}[{\quad\rm(a)}]
\item $\gdim_{\vf}(X)<\infty$.
\item $\gdim_{R'}(\comp X)<\infty$.
\item $\comp{X}$ is in $\A(R')$.
\item $\comp{X}$ is in $\A(\comp{R})$.
\end{enumerate}
When $R$ possesses a dualizing complex, these conditions are equivalent to:
\begin{enumerate}[{\quad\rm(e)}]
\item $X$ is in $\A(R)$.
\end{enumerate}
\end{prop}
\begin{proof}
  Indeed, (a) $\iff$ (b) by definition, while (b) $\iff$ (c)
  by~\cite[(3.1.10)]{christensen:gd}. Moving on, (c) $\iff$ (d) is contained
  in~\cite[(3.7.b)]{avramov:rhafgd}, and, when $R$ has a dualizing complex, the
  equivalence of (d) and (e) is~\cite[(3.7.a)]{avramov:rhafgd}.
\end{proof}


Now we turn to the behavior of G-dimension with respect to localizations.  Recall that,
given a prime ideal $\p$ and a totally reflexive $R$-module $G$, the $R_{\p}$-module
$G_{\p}$ is totally reflexive. From this it is clear that for any homologically finite
complex $W$, one has $\gdim_{R_{\p}}(W_{\p})\leq\gdim_R(W)$; see
\cite[(2.3.11)]{christensen:gd}.  For G-dimensions over $\vf$, we know only the following
weaker result;  see also~\cite[(10.2)]{avramov:homolhattfe}.  
Its proof is omitted for it is verbatim that of
\cite[(4.5)]{avramov:rhafgd}, which is the special case $X=S$; only, one uses
\ref{prop:likeAF4.3} instead of \cite[(4.3)]{avramov:rhafgd}.

\begin{prop}
\label{prop:localization}
Let $\vf\colon R\to S$ be a local homomorphism, $X$ a homologically finite complex of
$S$-modules. Let $\q$ be a prime ideal of $S$ and $\vf_\q$ the local homomorphism
$R_{\q\cap R}\to S_\q$.

If $\gdim_{\vf}(X)<\infty$, then $\gdim_{\vf_\q}(X_\q)<\infty$ under each of the
conditions:
\begin{enumerate}[{\quad\rm(1)}]
\item $\vf$ is essentially of finite type; or
\item $R$ has Gorenstein formal fibres. \qed
\end{enumerate}
\end{prop}

The next result shows that $\gdim_{\vf}(X)$ can be computed via any \textit{Gorenstein}
factorization of $\vf$, when such a factorization exists; see Definition~\ref{def:factor}.
When the factorization in the statement is regular, the equation becomes
$\gdim_{\vf}(X)=\gdim_{R'}(X)-\embdim(\Dot{\vf})$; compare this with
Definition~\ref{def:gdim}.

\begin{prop} \label{prop:comp-via-Gor-hom}
  Let $\vf\colon R\to S$ be a local homomorphism and $X$ a homologically finite complex of
  $S$-modules.  If $\vf$ possesses a Gorenstein factorization $R\xrightarrow{\Dot{\vf}}
  R'\xrightarrow{\vf'} S$, then
\[ \gdim_{\vf}(X)
=\gdim_{R'}(X)-\depth(\Dot{\vf}). \]
\end{prop}
\begin{proof}
  One may assume that $\HH(X)\neq 0$.  It is straightforward to verify that the diagram
  $R\to \comp{R'}\to\comp{S}$ is a Gorenstein factorization.  It follows
  from~\cite[(3.7)]{avramov:rhafgd} that $\comp{X}$ is in $\A(\comp{R})$ exactly when
  $\comp{X}$ is in $\A(\comp{R'})$, and, by Proposition \ref{prop:likeAF4.3}, this implies
  that $\gdim_{\vf}(X)$ is finite exactly when $\gdim_{R'}(X)$ is finite.  So one may
  assume that both the numbers in question are finite.  The Auslander-Bridger
  formula~\ref{thm:AB}, and the fact that $\depth_S(X)=\depth_{R'}(X)$, give the first of
  the following equalities:
\begin{align*}
  \gdim_{R'}(X)&=\gdim_{\vf}(X)+ [\depth(R')-\depth(R)] \\
  & = \gdim_{\vf}(X)+\depth(\Dot{\vf}).
\end{align*}
The second equality is by definition.
\end{proof}

\section{Projective dimension} \label{sec:pd}

In this section we introduce a new invariant: projective dimension over a local
homomorphism.  To begin with, one has the following proposition.  Its proof is similar to
that of Theorem \ref{thm:likeAF4.3}, and hence it is omitted.

\begin{prop}  \label{prop:def-pdim}
  Let $\vf\colon R\to S$ be a local homomorphism and $X$ a homologically finite complex of
  $S$-modules.  If $R\stackrel{\Dot{\vf_1}}{\to} R_1\stackrel{\vf_1'}{\to} \comp{S}$ and
  $R\stackrel{\Dot{\vf_2}}{\to} R_2\stackrel{\vf_2'}{\to} \comp{S}$ are Cohen
  factorizations of $\grave{\vf}$, then
\begin{xxalignat}{3}
  &{\hphantom{\square}}& \pdim_{R_1}(\comp{X})-\embdim(\Dot{\vf_1}) & =
  \pdim_{R_2}(\comp{X})-\embdim(\Dot{\vf_2}).  &&\square
\end{xxalignat}
\end{prop}

This leads to the following:

\begin{defn}  
  Let $\vf\colon R\to S$ be a local homomorphism and $X$ a homologically finite complex of
  $S$-modules.  The \textit{projective dimension of $X$ over $\vf$} is the quantity
\[ \pdim_{\vf}(X):= 
\pdim_{R'}(\comp{X})-\embdim(\Dot{\vf}) \] for some Cohen factorization $R\to R'\to
\comp{S}$ of $\grave{\vf}$.  The \textit{projective dimension of $\vf$} is defined to be
\[ \pdim(\vf):=\pdim_{\vf}(S).\]
\end{defn}

The first remark concerning this invariant is that there is an ``Auslander-Buchsbaum
formula'', which can be verified along the lines of its G-dimension counterpart, Theorem
\ref{thm:AB}.

\begin{property} \label{prop:pdAB}
  If $\pdim_{\vf}(X)<\infty$, then
\[ 
\pdim_{\vf}(X)=\depth(R)-\depth_S(X)
 \] 
\end{property}

Other basic rules that govern the behavior of this invariant can be read
from~\cite{avramov:solh}, although it was not defined there explicitly.  For instance,
\cite[(3.2)]{avramov:solh}, rather, its extension to complexes, see
\cite[(2.5)]{avramov:homolhattfe}, translates to

\begin{property}
\label{prop:pdvsfd}
There are inequalities:
\[
\fdim_R(X) - \embdim\vf \leq \pdim_{\vf}(X) \leq \fdim_R(X)\,.
\]
In particular, the finiteness of $\pdim_{\vf}(X)$ is independent of $S$ and $\vf$.
\end{property}

One can interpret the difference between $\fdim_R(X)$ and $\pdim_{\vf}(X)$ in terms of
appropriate depths:

\begin{prop}
 \label{prop:pd-fd}
 Let $\vf\colon R\to S$ be a local homomorphism and $X$ a homologically finite complex of
 $S$-modules.  Then
\[ 
\pdim_{\vf}(X) = \fd_R(X)+ \depth_R(X)-\depth_S(X).
\]
\end{prop}
\begin{proof}
  Indeed, by Property \ref{prop:pdvsfd}, we may assume that both $\fdim_R(X)$ and
  $\pdim_\vf(X)$ are finite.  Now, the first equality below is given
  by~\cite[(5.5)]{avramov:hdouc}, and the second is due to~\cite[(2.1)]{iyengar:dfcait}.
\[ 
\fd_R(X)=\sup(X\lotimes_R k)=\depth(R)-\depth_R(X).
\]
The Auslander-Buchsbaum formula~\ref{prop:pdAB} gives the desired formula.
\end{proof}

The G-dimension of a finite module, or a complex, is bounded above by its projective
dimension.  The same behavior carries over to modules and complexes over $\vf$.

\begin{prop} \label{prop:pd-gd-ineq}
  Let $\vf\colon R\to S$ be a local homomorphism.  For each homologically finite complex
  $X$ of $S$-modules, one has
\[ \gdim_{\vf}(X)\leq\pdim_{\vf}(X); \]
equality holds when $\pdim_{\vf}(X)<\infty$.
\end{prop}

\begin{proof}
  Let $R\xrightarrow{\Dot{\vf}} R'\to \comp{S}$ be a Cohen factorization of $\grave{\vf}$.
  Then
\[ \gdim_{\vf}(X)=\gdim_{R'}(\comp{X})-\embdim(\Dot{\vf})\leq 
\pdim_{R'}(\comp{X})-\embdim(\Dot{\vf}) =\pdim_{\vf}(X) \] with equality if
$\pdim_{R'}(\comp{X})$ is finite; see \cite[(2.3.10)]{christensen:gd}.
\end{proof}

Further results concerning $\pdim_{\vf}(\text{-})$ are given toward the end of the next
section.  One can introduce also Betti numbers and Poincar\'e series over local
homomorphisms; an in-depth analysis of these and related invariants is carried out
in~\cite{avramov:homolhattfe}.

\section{Ascent and descent of G-dimension} \label{sec:descent}

The heart of this section, and indeed of this paper, is the following theorem.  It is a
vast generalization of a stability result of Yassemi~\cite[(2.15)]{yassemi:gd}, and
contains Theorem B from the introduction.

\begin{thm} \label{thm:stability-tensor}
  Let $\vf\colon R\to S$ and $\sigma\colon S\to T$ be local homomorphisms.  Let $P$ be a
  complex of $T$-modules that is homologically finite with $\pdim_{\sigma}(P)$ finite.
  For every homologically finite complex $X$ of $S$-modules
\[ \gdim_{\sigma\vf}(X\lotimes_S P)=\gdim_{\vf}(X)+\pdim_{\sigma}(P).\]
In particular, $\gdim_{\sigma\vf}(X\lotimes_S P)$ and $\gdim_{\vf}(X)$ are simultaneously
finite.
\end{thm}

The theorem is proved in \ref{pf:stability-tensor}, toward the end of the section.  It is
worth remarking that the displayed formula is \emph{not} an immediate consequence of the
finiteness of the G-dimensions in question and appropriate Auslander-Bridger formulas.
What is missing is an extension of the Auslander-Buchsbaum formula \ref{prop:pdAB};
namely, under the hypotheses of the theorem above
\[
\pdim_{\sigma}(P)=\depth_S(X)-\depth_T(X\lotimes_SP)\,.
\]
It is not hard to deduce this equality from \cite[(2.4)]{foxby:daafuc}, using Cohen
factorizations; see the argument in \ref{pf:stability-tensor}.

We draw a few corollaries that illustrate the power of Theorem \ref{thm:stability-tensor}.
The first one is just the special case $X=S$ and $P=T$.

\begin{thm} \label{thm:gdim-compose}
  Let $\vf\colon R\to S$ and $\sigma\colon S\to T$ be local homomorphisms with
  $\pdim(\sigma)$ finite.  Then
\[ \gdim(\sigma\vf)=\gdim(\vf)+\pdim(\sigma).\]
In particular, $\gdim(\sigma\vf)$ is finite if and only if $\gdim(\vf)$ is finite.\qed
\end{thm}

The following example illustrates that the hypothesis on $\sigma$ cannot be weakened to
``$\gdim(\sigma)$ finite''.  A similar example is constructed
in~\cite[p.~931]{apassov:afm}.

\begin{ex} \label{ex:not-weaken}
  Let $R$ be a local, Cohen-Macaulay ring with canonical module $\omega$.  Set
  $S=R\ltimes\omega$, the ``idealization'' of $\omega$, and $\vf\colon R\to S$ the
  canonical inclusion.  Let $T=S/\omega\cong R$ with $\sigma\colon S\to T$ the natural
  surjection.
  
  Now, $\sigma\vf=\mathrm{id}^R$, hence $\gdim(\sigma\vf)=0$, for example, by Proposition
  \ref{prop:pd-gd-ineq}; also, $S$ is Gorenstein~\cite[(3.3.6)]{bruns:cmr}, so
  $\gdim(\sigma)$ is finite.
  
  We claim that $\gdim(\vf)$ is finite if and only if $R$ is Gorenstein.  Indeed,
  $\gdim(\vf)$ and $\gdim_R(S)$ are simultaneously finite, by Corollary
  \ref{cor:vf-is-finite}.  From~\ref{subprop:3} it follows that $\gdim_R(S)<\infty$ if and
  only if $\gdim_R(\omega)<\infty$.  The finiteness of $\gdim_R(\omega)$ is equivalent to
  $R$ being Gorenstein~\cite[(3.4.12)]{christensen:gd}.
\end{ex}

As noted in the introduction, Theorem \ref{thm:gdim-compose} allows one to extend certain
results of Avramov and Foxby on (quasi-)Gorenstein homomorphisms.  In order to describe
these, and because they are required in the sequel, we recall the relevant notions.

\begin{para} 
\label{para:quasi-Gorenstein}
\label{para:1}
Let $R$ be a local ring with residue field $k$.  The \textit{Bass series of $R$} is the
formal power series $I_R(t)=\sum_i\mu^i_R(R)t^i$ where $\mu_R^i(R)=\rank_k\ext_R^i(k,R)$.
An important property of the Bass series is that $R$ is Gorenstein if and only if $I_R(t)$
is a polynomial~\cite[(18.1)]{matsumura:crt}.
  
Let $\vf\colon (R,\m)\to (S,\n)$ be a local homomorphism of finite G-dimension.  Let
$I_{\vf}(t)$ denote the \textit{Bass series of $\vf$}, introduced in ~\cite[Section
7]{avramov:rhafgd}. The Bass series is a formal Laurent series with nonnegative integer
coefficients and satisfies the equality
\begin{equation}  \label{eqn:bass}
I_R(t)I_{\vf}(t)=I_S(t). \tag{$\dagger$}
\end{equation}
Let $\sigma\colon (S,\n)\to (T,\p)$ also be a homomorphism of finite G-dimension.
Assuming $\gdim(\sigma\vf)$ is finite as well, it follows from (\ref{eqn:bass}) that
\begin{equation} \label{eqn:more-bass}
I_{\sigma\vf}(t)=I_{\sigma}(t)I_{\vf}(t). \tag{$\ddagger$}
\end{equation}

The homomorphism $\vf$ is said to be \textit{quasi-Gorenstein at $\n$} if $I_{\vf}(t)$ is
a Laurent \textit{polynomial}.  When $\pdim(\vf)<\infty$ and $\vf$ is quasi-Gorenstein at
$\n$, one says that $\vf$ is \emph{Gorenstein at $\n$};
see~\cite[(7.7.1)]{avramov:rhafgd}.

A noteworthy aspect of the class of such homomorphisms is that it is closed under
compositions: if $\vf$ and $\sigma$ are quasi-Gorenstein at $\n$ and $\p$, respectively,
then $\sigma\vf$ is quasi-Gorenstein at $\p$.  This follows
from~\cite[(7.10)]{avramov:rhafgd} and (\ref{eqn:more-bass}) above.
\end{para}

The result below is a decomposition theorem for Gorenstein and quasi-Gorenstein
homomorphisms; it is Theorem C announced in the introduction.

\begin{thm} \label{thm:fin-fd-descent}
  Let $\vf\colon (R,\m)\to (S,\n)$ and $\sigma\colon (S,\n)\to (T,\p)$ be local
  homomorphisms with $\pdim(\sigma)$ finite.  If $\sigma\vf$ is (quasi-)Gorenstein at
  $\p$, then $\vf$ is (quasi-)Gorenstein at $\n$ and $\sigma$ is Gorenstein at $\p$.
\end{thm}

\begin{proof}
  Assume that $\sigma\vf$ is quasi-Gorenstein at $\p$; so, it has finite
  G-dimension, and $I_{\sigma\vf}(t)$ is a Laurent polynomial.  Now,
  $\gdim(\vf)$ is finite, by Theorem~\ref{thm:gdim-compose}, as is $\gdim(\sigma)$, by
  hypothesis, so equality (\ref{eqn:more-bass}) in \ref{para:quasi-Gorenstein} applies to yield an
  equality of formal Laurent series
\[ 
I_{\sigma\vf}(t)=I_{\sigma}(t)I_{\vf}(t)
 \]
 In particular, $I_{\sigma}(t)$ and $I_{\vf}(t)$ are Laurent polynomials as well. Thus,
 both $\sigma$ and $\vf$ are quasi-Gorenstein at the appropriate maximal ideals.
 Moreover, $\sigma$ is Gorenstein because $\pdim(\sigma)$ is finite.
 
 Suppose that $\sigma\vf$ is Gorenstein at $\p$.  Since $\pdim(\sigma\vf)$ and
 $\pdim(\sigma)$ are both finite, \cite[(3.2)]{foxby:daafuc}, in conjunction with
 Proposition~\ref{prop:pd-fd}, yields that $\pdim(\vf)$ is finite. The already established
 part of the theorem gives the desired conclusion.
\end{proof}

The next theorem generalizes another stability result of
Yassemi~\cite[(2.14)]{yassemi:gd}.

\begin{thm} \label{thm:stability-hom}
  Let $\vf\colon R\to S$ be a local homomorphism and $P$ a homologically finite complex of
  $S$-modules with $\pdim_S(P)$ finite.  For every homologically finite complex $X$ of
  $S$-modules
\begin{align*}
  \gdim_{\vf}(\rhom_S(P,X))&=\gdim_{\vf}(X)-\inf(P).
\end{align*}
Thus, $\gdim_{\vf}(X)$ and $\gdim_{\vf}(\rhom_S(P,X))$ are simultaneously finite.
\end{thm}

\begin{proof}
  The tensor evaluation morphism $X\lotimes_S\rhom_S(P,S)\to\rhom_S(P,X)$ is an
  isomorphism, as $P$ has finite projective dimension.  Thus
\begin{align*}
  \gdim_{\vf}(\rhom_S(P,X))&=\gdim_{\vf}(X\lotimes_S\rhom_S(P,S))\\
  &=\gdim_{\vf}(X)+\pdim_S(\rhom_S(P,S))\\
  &=\gdim_{\vf}(X)-\inf(P)
\end{align*}
where the second equality follows from Theorem~\ref{thm:stability-tensor} because
$\rhom_S(P,S)$ has finite projective dimension over $S$.
\end{proof}

Next we record the analogue of Theorem~\ref{thm:stability-tensor} for projective
dimension; its proof is postponed to \ref{pf:stab-pd-hom}.

\begin{thm} \label{thm:stab-pd-tensor}
  Let $\vf\colon R\to S$ and $\sigma\colon S\to T$ be local homomorphisms.  Let $P$ be a
  complex of $T$-modules that is homologically finite with $\pdim_{\sigma}(P)$ finite.
  For every homologically finite complex $X$ of $S$-modules
\[ \pdim_{\sigma\vf}(X\lotimes_S P)=\pdim_{\vf}(X)+\pdim_{\sigma}(P).\]
In particular, $\pdim_{\sigma\vf}(X\lotimes_S P)$ and $\pdim_{\vf}(X)$ are simultaneously
finite.
\end{thm}

Finally, here is the analogue of Theorem~\ref{thm:stability-hom}; it can be deduced
from~\ref{thm:stab-pd-tensor} in the same way that~\ref{thm:stability-hom} was deduced
from~\ref{thm:stability-tensor}.

\begin{thm} \label{thm:stab-pd-hom}
  Let $\vf\colon R\to S$ be a local homomorphism and $P$ a homologically finite complex of
  $S$-modules with $\pdim_S(P)$ finite.  For every homologically finite complex $X$ of
  $S$-modules
\begin{align*}
  \pdim_{\vf}(\rhom_S(P,X))&=\pdim_{\vf}(X)-\inf(P).
\end{align*}
In particular, $\pdim_{\vf}(X)$ and $\pdim_{\vf}(\rhom_S(P,X))$ are simultaneously
finite.\qed
\end{thm}

The proof of Theorem \ref{thm:stability-tensor} uses a convenient construction,
essentially given in \cite{avramov:solh}, of Cohen factorizations of compositions of local
homomorphisms.

\begin{para} \label{para:big-diagram}
  Let $R\xrightarrow{\vf} S\xrightarrow{\sigma} T$ be local homomorphisms, and let
\[ 
R\xrightarrow{\Dot{\vf}} R'\xrightarrow{\vf'}S 
  \qquad\text{and}\qquad 
       R'\xrightarrow{\Dot{\rho}}R''\xrightarrow{\rho'}T 
\] 
be regular factorizations of $\vf$ and $\sigma\vf'$, respectively.  The map $\rho'$
factors through the tensor product $S'=R''\otimes_{R'}S$ giving the following commutative
diagram
\[
\xymatrixrowsep{2.5pc}
\xymatrixcolsep{2.5pc} 
\xymatrix { 
  && R'' \ar@{->}[dr]^{\vf''}  \ar@/^2.5pc/[ddrr]^{\rho'=\sigma'\vf''}  \\
  & R' \ar@{->}[dr]^{\vf'} \ar@{->}[ur]^{\dot\rho} && S' \ar@{->}[dr]^{\sigma'}  \\
  R \ar@{->}[rr]^{\vf} \ar@/^2.5pc/[uurr]^{\dot\rho \dot\vf} \ar@{->}[ur]^{\dot\vf} && S
  \ar@{->}[rr]^{\sigma} \ar@{->}[ur]^{\dot\sigma} && T }
\]
where $\dot\sigma$ and $\vf''$ are the natural maps to the tensor products.  Then the
diagrams $S\to S'\to T$, $R'\to R''\to S'$, and $R\to R''\to T$ are regular factorizations
with
\[
\embdim(\Dot{\sigma})=\embdim(\Dot{\rho}) \qquad\text{and}\qquad
\embdim(\Dot{\rho}\Dot{\vf}) =\embdim(\Dot{\rho})+\embdim(\Dot{\vf}).
\]
Indeed, by flat base change, $\Dot{\sigma}$ is flat and has closed fiber $S'\otimes_S
l=R''\otimes_{R'} l$, which is regular.  Here $l$ is the common residue field of $R'$ and
$S$.  This tells us that $S\to S'\to T$ is a regular factorization and that
$\embdim(\Dot{\sigma})=\embdim(\Dot{\rho})$.

The diagram $R'\to R''\to S'$ is a regular factorization because $\Dot{\rho}$ is flat with
a regular closed fibre, by hypothesis, and $\vf''$ is surjective, by base change.

As to the diagram $R\to R''\to T$, let $\m$ and $\m'$ denote the maximal ideals of $R$ and $R'$,
respectively.  The induced map $R'/\m R'\to R''/\m R''$ is flat with closed fibre
$R''/\m'R''$.  Since $R'/\m R'$ and $R''/\m'R''$ are both regular, the same is true of
$R''/\m R''$, by~\cite[(2.2.12)]{bruns:cmr}.  Thus, $R\to R''\to T$ is a regular
factorization.  Furthermore, it is stated explicitly in the proof of \emph{loc.~cit.} that
$\embdim(R''/\m R'')=\embdim(R'/\m R')+\embdim(R''/\m' R'')$, which explains the second
formula above.
\end{para}

\begin{para}\textit{Proof of Theorem~\ref{thm:stability-tensor}.}
\label{pf:stability-tensor}
Note that $X\lotimes_S P$ is homologically finite over $T$ by Lemma~\ref{lem:back-finite},
so one may speak of its G-dimension over $\sigma\vf$.  Passing to the completions of $S$
and $T$ at their respective maximal ideals, and replacing $X$ and $P$ by
$\comp{S}\otimes_S X$ and $\comp{T}\otimes_T P$, respectively, one may assume that $S$ and
$T$ are complete.  In doing so, one uses the isomorphism
\[ (\comp{S}\otimes_S X)\lotimes_{\comp{S}} (\comp{T}\otimes_T P)
\simeq \comp{T}\otimes_T (X\lotimes_S P). \]

The next step is the reduction to the case where $\vf$ and $\sigma$ are surjective.  To
achieve this, take Cohen factorizations $R\to R'\to S$ and $R'\to R''\to T$, and expand to
a commutative diagram as in~\ref{para:big-diagram}.

Let $X'=S'\otimes_S X$.  Since $S'=R''\otimes_{R'}S$, by construction, $X'\cong
R''\otimes_{R'} X$ and hence $X'\lotimes_{S'}P\simeq X\lotimes_S P$.  Since $R'\to R''$ is
faithfully flat, \ref{lem:fflat} yields
\[ \gdim_{R'}(X)=\gdim_{R''}(X'). \]
The preceding equality, in conjunction with those in~\ref{para:big-diagram}, yields
\begin{align*}
  \pdim_{\sigma}(P) & = \pdim_{S'}(P)-\embdim(\Dot{\rho}) \\
  \gdim_{\vf}(X) & =\gdim_{R''}(X')-\embdim(\Dot{\vf}) \\
  \gdim_{\sigma\vf}(X\lotimes_S P) & = \gdim_{R''}(X'\lotimes_{S'}P)-\embdim(\Dot{\rho})
  -\embdim(\Dot{\vf})
\end{align*}
Therefore, it suffices to verify the identity for the diagram $R''\to S'\to T$ and
complexes $X'$ and $P$.  This places us in the situation where $R\to S$ is surjective, $P$
is homologically finite over $R$, and then the equality we seek is
\[ \gdim_R(X\lotimes_S P)=\gdim_R(X)+\pdim_S(P).\]

It suffices to prove that the G-dimensions over $R$ of $X$ and of $X\lotimes_S P$ are
simultaneously finite.  For, when they are both finite, one has
\begin{align*}
  \gdim_R(X\lotimes_S P)
  & = \depth(R)-\depth_S(X\lotimes_S P)  \\
  & = \depth(R)-\depth_S(X)+\pdim_S(P) \\
  & = \gdim_R(X)+\pdim_S(P)
\end{align*}
where the first and the third equalities are by the Auslander-Bridger formula, while the
one in the middle is by~\cite[(2.2)]{iyengar:dfcait}.

The rest of the proof is dedicated to proving that $X$ and $X\lotimes_S P$ have finite
G-dimension over $R$ simultaneously.  In view of~\ref{para:reflexive-2}, this is
tantamount to proving:
\begin{enumerate}[{\quad\rm(a)}]
\item $\rhom_R(X,R)$ is homologically bounded if and only if the same is true of
  $\rhom_R(X\lotimes_S P,R)$; and
\item the biduality morphisms $\delta^R_X$ and $\delta^R_{X\lotimes_S P}$, defined as
  in~\ref{para:reflexive}, are isomorphisms simultaneously.
\end{enumerate}

The proofs of (a) and (b) use the following observation: when $U$ and $V$ are complexes of
$S$-modules such that $V$ is homologically finite and $\pdim_S(V)<\infty$, the natural
morphism
\begin{equation} \label{eqn:1}
\theta_{UV}\colon\rhom_R(U,R)\lotimes_S\rhom_S(V,S)\to 
   \rhom_R(U\lotimes_S V,R) \tag{$\ast$}
\end{equation}
is an isomorphism.  Indeed, it is the composition of tensor evaluation
\[ 
\rhom_R(U,R)\lotimes_S\rhom_S(V,S) \to \rhom_S(V, \rhom_R(U, R)), 
 \] 
which is an isomorphism for $V$ as above, followed by
adjunction
\[ 
\rhom_S(V, \rhom_R(U, R))\xrightarrow{\simeq} \rhom_R(U\lotimes_S V,R). 
\]

\textit{Proof of (a)}.  Since $\rhom_S(P,S)$ has finite projective dimension over $S$, one
has the isomorphism
\[ 
\te_{XP}\colon\rhom_R(X,R)\lotimes_S\rhom_S(P,S)\to
\rhom_R(X\lotimes_S P,R). 
 \] Thus, Theorem~\ref{thm:GAI} implies the desired equivalence.

\textit{Proof of (b)}.  Consider the following commutative diagram of morphisms of
complexes of $S$-modules.
\[ 
\xymatrix{
  X\lotimes_S P \ar[r]^(0.3){(\delta_X^R)\lotimes_S P} \ar@{=}[dd] &
  \rhom_R(\rhom_R(X,R),R)\lotimes_S P \ar[d]_{\simeq}^{\nu} \\
  & \rhom_R(\rhom_R(X,R)\lotimes_S\rhom_S(P,S),R) \\
  X\lotimes_S P \ar[r]^(0.3){\delta^R_{(X\lotimes_S P)}} & \rhom_R(\rhom_R(X\lotimes_S
  P,R),R) \ar[u]^{\simeq}_{\rhom_R(\theta_{X\!P},R)} } 
\] 
The morphism $\nu$ is the composition $\te_{UV}\circ (1\lotimes_S\delta^S_P)$ where
$U=\rhom_R(X,R)$ and $V=\rhom_S(P,S)$.  Note that $\delta^S_P$, and hence
$1\lotimes_S\delta^S_P$, is an isomorphism because $\pdim_S(P)$ is finite.  Furthermore,
$\te_{UV}$ is an isomorphism since $\pdim_S(V)$ is finite.  This is why $\nu$ is an
isomorphism.

From the diagram one obtains that $\delta_{X\lotimes_S P}^R$ is an isomorphism if and only
if $(\delta^R_X)\lotimes_S P$ is.  By Proposition~\ref{prop:tool1}, the morphisms
$(\delta^R_X)\lotimes_S P$ and $\delta^R_X$ are isomorphisms simultaneously, as
$\pdim_S(P)$ is finite.  \qed
\end{para}

To wrap up this section, we give the

\begin{para}\textit{Proof of Theorem \ref{thm:stab-pd-tensor}.}
\label{pf:stab-pd-hom}
Arguing as in the proof of Theorem~\ref{thm:stability-tensor} one reduces to the case
where $\vf$ and $\sigma$ are surjective and $\pdim_S(P)$ is finite.  In this situation,
one has to verify that $\pdim_{R}(X)$ and $\pdim_{R}(X\lotimes_S P)$ are simultaneously
finite.  Let $k$ be the residue field of $R$.  It suffices to show that $\amp(k\lotimes_R
X)$ and $\amp(k\lotimes_R(X\lotimes_S P))$ are simultaneously finite.  By the isomorphism
 \[
 k\lotimes_R(X\lotimes_S P)\simeq(k\lotimes_R X)\lotimes_S P\,,
\]
this follows from Theorem~\ref{thm:GAI}. \qed
\end{para}

\section{Detecting the Gorenstein property} \label{sec:apps}

The theorem below extends the Auslander-Bridger
characterization~\cite[(4.20)]{auslander:smt} of Gorenstein rings.

\begin{thm} \label{thm:gor-char-1}
  Let $R$ be a local ring. The following conditions are equivalent.
\begin{enumerate}[{\quad\rm(a)}] 
\item $R$ is Gorenstein.
\item For every local homomorphism $\vf\colon R\to S$ and for every homologically finite
  complex $X$ of $S$-modules, $\gdim_{\vf}(X)<\infty$.
\item There is a local homomorphism $\vf\colon R\to S$ and an ideal $I$ of $S$ such that
  $I\supseteq \m S$, where $\m$ is the maximal ideal of $R$, and
  $\gdim_{\vf}(S/I)<\infty$.
\end{enumerate}
\end{thm}

\begin{proof}
  ``(a)$\implies$(b)''.  Let $R\to R'\to \comp{S}$ be a Cohen factorization of
  $\grave{\vf}$.  The $R'$-module $\HH(\comp{X})$ is finite, because the $S$-module
  $\HH(X)$ is finite.  Since $R$ is Gorenstein, so is $R'$~\cite[(3.3.15)]{bruns:cmr}.
  Thus, $\gdim_{R'}(\comp{X})<\infty$, that is to say, $\gdim_{\vf}(X)<\infty$; see
  Proposition~\ref{prop:likeAF4.3}.
  
  ``(b)$\implies$(c)'' is trivial.
  
  ``(c)$\implies$(a)''.  Let $R\to R'\to \comp{S}$ be a Cohen factorization.  Composing
  with the surjection $\comp{S}\xrightarrow{\pi}\comp{S}/I\comp{S}$ gives a diagram $R\to
  R'\to \comp{S}/I\comp{S}$ that is also a Cohen factorization.  Since
  $\gdim_{R'}(\comp{S}/I\comp{S})$ is finite, so is $\gdim(\pi\grave{\vf})$.  The
  composition $\pi\grave{\vf}$ factors through the residue field $k$ of $R$, giving the
  commutative diagram:
\[ 
\xymatrix{
  R \ar[dr] \ar[rr] & & \comp{S}/I\comp{S} \\
  & k \ar[ur] } 
\] 
The map $k\to \comp{S}/I\comp{S}$ has finite projective dimension because $k$ is a field.
Therefore, Theorem~\ref{thm:gdim-compose} implies that the surjection $R\to k$ has finite
G-dimension.  Thus, $R$ is Gorenstein by~\ref{subprop:4} and~\cite[(4.20)]{auslander:smt}.
\end{proof}

When $\vf$ is finite and $X$ is a module of finite projective dimension over both $R$ and
$S$, the implication ``(c)$\implies$(a)'' in the next result was proved by
Apassov~~\cite[Theorem G']{apassov:afm}.

\begin{thm} \label{thm:like-apassov}
  Let $\vf\colon R\to S$ be a local homomorphism such that $S$ is Gorenstein.  The
  following conditions are equivalent.
\begin{enumerate}[{\quad\rm(a)}] 
\item $R$ is Gorenstein.
\item $\gdim(\vf)$ is finite.
\item There exists a homologically finite complex $P$ of $S$-modules such that
  $\pdim_S(P)$ is finite and $\gdim_{\vf}(P)$ is finite.
\end{enumerate}
\end{thm}

\begin{proof} The implication
  ``(a)$\implies$(b)'' is contained in Theorem~\ref{thm:gor-char-1}, while
  ``(b)$\iff$(c)'' is given by Theorem~\ref{thm:stability-tensor}.
  
  ``(b)$\implies$(a)''. As $S$ is Gorenstein, $I_S(t)=t^{\dim(S)}$, and hence
\[
I_R(t) I_{\vf}(t) = I_S(t)=t^{\dim(S)}\,,
 \]
 by equality (\ref{eqn:bass}) in~\ref{para:1}.  Now, both $I_R(t)$ and $I_{\vf}(t)$ are
 Laurent series with nonnegative coefficients, so that $I_R(t)$ is a polynomial. This, as
 noted in~\ref{para:1}, implies that $R$ is Gorenstein.
\end{proof}

The last theorem in this section is a characterization of the Gorenstein property of a
local ring in terms of the finiteness of G-dimension of Frobenius-like endomorphisms.  In
order to describe this, we recall the definition of an invariant introduced by Koh and
Lee~\cite[(1.1)]{koh:nionlr}.

\begin{para} \label{para:col(R)}
For a finite module $M$ over a local ring $(S,\n)$, set
\[ 
\s(M)=\inf\{t\geq 1\mid\soc(M)\not\subseteq\n^t M\}
\]
where $\soc(M)$ is the socle of $M$.  Furthermore, let
\[ 
\crs(S)=\inf\{\s(S/(\x))\mid 
  \text{$\x=x_1,\ldots,x_r$ is a maximal $S$-sequence}\}.  
\]
\end{para}

The following is a complex version of Koh-Lee~\cite[(2.6)]{koh:srotmimr} (see also
Miller~\cite[(2.2.8)]{miller:tfeahd}), which, in turn, generalizes a theorem of
J.~Herzog~\cite[(3.1)]{herzog:rdcpuf}.

\begin{prop} \label{prop:koh-lee}
  Let $\vf\colon (R,\m)\to (S,\n)$ be a local homomorphism for which
  $\vf(\m)\subseteq\n^{\crs(S)}$, and $X$ a homologically finite complex of $R$-modules.
  If there is an integer $t\geq\sup(X)$ such that $\tor_{t+i}^R(X,S)=0$ for $1\leq
  i\leq\depth(S)+2$, then $\pdim_R(X)<\infty$.
\end{prop}

\begin{proof} 
  Replace $X$ with a minimal $R$-free resolution to assume that each $X_i$ is a finite
  free $R$-module, and $\partial(X)\subseteq \m X$.  Set $Y=X\otimes_R S$.  Then $Y$ is a
  complex of finite free $S$-modules with $\partial(Y)\subseteq\n^{\crs(S)}Y$ and
  $\HH_i(Y)=\tor^R_i(X,S)$.
  
  The desired conclusion is that $X_i=0$ for $i\gg 0$.  By the minimality of $X$, it
  suffices to prove that $X_i=0$, equivalently, $Y_i=0$, for some $i>\sup(X)$.  Let
  $r=\depth(S)$ and $C=\coker(\partial^Y_{t+1})$.  The truncated complex
\[ 
Y_{t+r+2}\xrightarrow{\partial_{t+r+2}} Y_{t+r+1}
  \xrightarrow{\partial_{t+r+1}}\cdots\xrightarrow{\partial_{t+1}} Y_t\to 0 
\]
is the beginning of a minimal $S$-free resolution of $C$.  If $\pdim_S(C)=\infty$,
then~\cite[(2.2.5),(2.2.6)]{miller:tfeahd} implies that each row of $\partial^Y_{t+r+2}$
has an entry outside $\n^{\crs(S)}$, a contradiction.  Thus, $\pdim_S(C)<\infty$, and the
Auslander-Buchsbaum formula implies that the $\pdim_S(C)\leq r$.  The minimality of the
complex above implies that $Y_{t+r+1}=0$, completing the proof.
\end{proof}

An arbitrary local homomorphism of finite G-dimension is far from being quasi-Gorenstein.
Indeed, when $R$ is Gorenstein, any local homomorphism $\vf\colon R\to S$ has finite
G-dimension, see Theorem \ref{thm:gor-char-1}, whereas, by \cite[(8.2)]{avramov:rhafgd},
such a $\vf$ is quasi-Gorenstein if and only if $S$ is Gorenstein. Endomorphisms however
are much better behaved in this regard.

\begin{prop} \label{prop:5-6}
  Let $\vf\colon (R,\m) \to (R,\m)$ be a local homomorphism.  If $\gdim(\vf)$ is finite,
  then $\vf^n$ is quasi-Gorenstein at $\m$, for each integer $n\geq1$.
\end{prop}

If the finiteness of G-dimension localizes--see Proposition \ref{prop:localization}--then
one could draw the stronger conclusion that $\vf^n$ is quasi-Gorenstein at each prime
ideal.

\begin{proof}[Proof of Proposition \ref{prop:5-6}]
  Suppose that $\gdim(\vf)$ is finite. The equality (\ref{eqn:bass}) in~\ref{para:1}
  yields $I_{\vf}(t)=1$ so that $\vf$ is quasi-Gorenstein at $\m$, the maximal ideal of
  $R$.  In the light of the discussion in~\ref{para:quasi-Gorenstein}, the same is true of
  the $n$-fold composition $\vf^{n}$, for all integers $n\geq 1$.
\end{proof}

We are now ready to prove the following theorem that subsumes Theorem A in the
introduction.

\begin{thm}  \label{thm:gdim(f)-finite}
  Let $\vf\colon (R,\m)\to (R,\m)$ be a local homomorphism such that
  $\vf^i(\m)\subseteq\m^2$ for some integer $i\geq 1$.  The following conditions are
  equivalent.
\begin{enumerate}[{\quad\rm(a)}] 
\item The ring $R$ is Gorenstein.
\item $\gdim(\vf^n)$ is finite for some integer $n\geq 1$.
\item There is a homologically finite complex $P$ of $R$-modules with $\pdim_R(P)$ finite
  and $\gdim_{\vf^n}(P)$ finite, for some integer $n\geq 1$.
\end{enumerate}
When these conditions hold, $\gdim(\vf^m)=0$, for all $m\geq 1$.
\end{thm}

\begin{proof}
  ``(a)$\implies$(b)'' is contained in Theorem~\ref{thm:gor-char-1}.
  
  ``(b)$\implies$(c)'' is trivial.
  
  ``(c)$\implies$(a)''.  By Theorem~\ref{thm:stability-tensor}, $\gdim(\vf^n)$ is finite.
  The completion $\comp{\vf}\colon \comp{R}\to\comp{R}$ is an endomorphism of $\comp{R}$
  such that $\comp{\vf}^i(\comp{\m})\subseteq\comp{\m}^2$ and
  $\comp{\vf^n}=(\comp{\vf})^n$.  Furthermore, $R$ is Gorenstein if and only if $\comp{R}$
  is Gorenstein, and by~\ref{subprop:1}, $\gdim(\vf^n)$ is finite if and only if
  $\gdim(\comp{\vf}^n)$ is finite.  Thus, passing to $\comp{R}$, one may assume that $R$
  is complete.  Hence, $R$ has a dualizing complex $D$.
  
  Since $\gdim(\vf^n)$ is finite, Proposition~\ref{prop:5-6} implies that the $s$-fold
  composition $\vf^{sn}$ of $\vf^n$ is also quasi-Gorenstein at $\m$ for all integers
  $s\geq 1$.  Thus, $D\lotimes_R {}^{\vf^{sn}}\!R$ is a dualizing complex for $R$, for
  each $s\geq 1$, by~\cite[(7.8)]{avramov:rhafgd}.  This implies that $\HH(D\lotimes_R
  {}^{\vf^{sn}}\!R)$ is finite; that is to say, $\tor_i^R(D,{}^{\vf^{sn}}\!R)=0$ for all
  $i\gg 0$.  Therefore, $\pdim_R(D)$ is finite, by Proposition~\ref{prop:koh-lee}.  This
  is equivalent to $R$ being Gorenstein; see, for example~\cite[(3.4.12)]{christensen:gd}.
  This completes the proof that (c) implies (a).
  
  When these conditions hold, the Auslander-Bridger formula~\ref{thm:AB} gives
\[ 
\gdim(\vf^m)=\depth(R)-\depth(R)=0.
  \]
This is the desired formula.
\end{proof}

The preceding theorem raises the problem: given a local ring $(R,\m)$ construct
endomorphisms of $R$ that map $\m$ into $\m^2$.  The prototype is the Frobenius
endomorphism of a local ring of characteristic $p$.  There are many such endomorphisms of
power series rings over fields.  The following example gives a larger class of complete
local rings with nontrivial endomorphisms.

\begin{ex} 
  Let $k$ be a field and $X_1,\ldots,X_n$ analytic indeterminates and $F_1,\ldots,F_m\in
  k[X_1,\ldots,X_n]$ homogeneous polynomials, and set
\[ 
R=k[\![X_1,\ldots,X_n]\!]/(F_1,\ldots,F_m)
=k[\![x_1,\ldots,x_n]\!].
\]
Let $g$ be an element in $(x_1,\ldots,x_n)R$.  The assignment $x_i\mapsto x_ig$ gives rise
to a well-defined ring endomorphism $\vf$ of $R$ such that $\vf(\m)\subseteq \m^2$.

One property of the Frobenius endomorphism that is hard to mimic is the finiteness of the
length of $R/\vf(\m)R$.  Again, over power series rings such endomorphisms abound.  The
desired property is satisfied by the ring $R$ constructed above, when it is a
one-dimensional domain and $g\ne0$.  Examples in dimension two or higher can be built from
these by considering $R[\![Y_1,\ldots,Y_m]\!]$.

More interesting endomorphisms can be obtained as follows.  Let
\[ R=k[\![X_1,\ldots,X_n]\!]/(G_1-H_1,\ldots,G_m-H_m) \]
where, for each $i$, the elements $G_i$ and $H_i$ are monomials of the same total degree.
For each positive integer $t$, the assignment $x_i\mapsto x_i^t$ gives rise to a ring
endomorphism $\vf_t$ of $R$ such that $\vf_t(\m)\subseteq\m^t$ and $R/\vf_t(\m)R$ has
finite length.  This method allows one to construct Cohen-Macaulay normal domains of
arbitrarily large dimension with nontrivial endomorphisms; consider, for example, the
maximal minors of a $2\times r$ matrix of variables.
\end{ex}

\section{Finiteness of G-dimension over $\vf$}
The import of the results of this section is that the finiteness of $\gdim_\vf(X)$ is
intrinsic to the $R$-module structure on $X$; this is exactly analogous to the behavior of
$\pdim_{\vf}(X)$; see \ref{prop:pdvsfd}. When $R$ is complete, it is contained in
Proposition \ref{prop:likeAF4.3}; see also Theorem \ref{thm:gdvsgfd} ahead.

\begin{thm} \label{thm:indep-of-vf}
  Let $\vf\colon R\to S$ and $\psi\colon R\to T$ be local homomorphisms.  Let $X$ and $Y$
  be homologically finite complexes of $S$-modules and $T$-modules, respectively, that are
  isomorphic in the derived category of $R$.  Then $\gdim_{\vf}(X)$ is finite if and only
  if $\gdim_{\psi}(Y)$ is finite.
\end{thm}

\begin{proof}
  One may assume that $X$ and $Y$ are homologically nonzero.  First, we reduce to the case
  where $\m$, the maximal ideal of $R$, annihilates $\HH(X)$ and $\HH(Y)$.  To this end,
  let $K$ be the Koszul complex on a finite generating sequence for $\m$.  Since
  $X\otimes_R K=X\otimes_S(S\otimes_R K)$ and $S\otimes_R K$ is a finite free complex of
  $S$-modules, Theorem~\ref{thm:stability-tensor} yields that $\gdim_{\vf}(X\otimes_R K)$
  and $\gdim_{\vf}(X)$ are simultaneously finite.  Similarly, $\gdim_{\psi}(Y\otimes_R K)$
  and $\gdim_{\psi}(Y)$ are simultaneously finite.  Moreover, $X\otimes_R K$ and
  $Y\otimes_R K$ are isomorphic in the derived category of $R$.  As $\m$ annihilates
  $\HH(X\otimes_R K)$ and $\HH(Y\otimes_R K)$--see, for instance,
  \cite[(1.2)]{iyengar:dfcait}--replacing $X$ and $Y$ with $X\otimes_R K$ and $Y\otimes_R
  K$, respectively, gives the desired reduction.
  
  Let $\alpha\colon X\to Y$ be an isomorphism.  Let $\wt{\vf}\colon\comp{R}\to\wt{S}$ and
  $\wt{\psi}\colon\comp{R}\to\wt{T}$ be the $\m$-adic completions of $\vf$ and $\psi$,
  respectively, and $\iota\colon R\to\comp{R}$ the completion map.  In the derived
  category of $R$, one has a commutative diagram:
\[ 
\xymatrix{ X=R\otimes_R X \ar[r]^(0.55){\iota\otimes_R 1}_(0.55){\simeq}
  \ar[d]^\alpha_{\simeq} & \comp{R}\otimes_R X \ar[r]^{\wt{\vf}\otimes_{\vf}1}_{\simeq}
  \ar[d]^{1\otimes_R \alpha}_{\simeq} & \wt{S}\otimes_S X \\
  Y=R\otimes_R Y \ar[r]^(0.55){\iota\otimes_R 1}_(0.55){\simeq} & \comp{R}\otimes_R Y
  \ar[r]^{\wt{\psi}\otimes_{\psi}1}_{\simeq} & \wt{T}\otimes_T Y }
\] 
Both $R\to\comp{R}$ and $S\to\wt{S}$ are flat, so at the level of homology the top row of
the diagram reads $\HH(X)\to\comp{R}\otimes_R \HH(X)\to\wt{S}\otimes_S\HH(X)$.  Since
$\HH(X)$ is annihilated by $\m$, these are both bijective, that is, $\iota\otimes_R 1$ and
$\wt{\vf}\otimes_{\vf} 1$ are isomorphisms.  A similar reasoning justifies the
isomorphisms in the bottom row.

In the biimplications below, the first is by~\ref{subprop:1}, the second is by
Proposition~\ref{prop:likeAF4.3}, while the third is due to the fact that, by the diagram
above, $\wt{S}\otimes_S X$ and $\comp{R}\otimes_R X$ are isomorphic.
\begin{align*}
  \text{$\gdim_{\vf}(X)$ is finite}
  &\iff \text{$\gdim_{\wt{\vf}}(\wt{S}\otimes_S X)$ is finite} \\
  &\iff \text{$\wt{S}\otimes_S X$ is in $\A(\comp{R})$} \\
  &\iff \text{$\comp{R}\otimes_R X$ is in $\A(\comp{R})$}
\end{align*}
By the same token, $\gdim_{\psi}(Y)$ is finite if and only if $\comp{R}\otimes_R Y$ is in
$\A(\comp{R})$.  This gives the desired conclusion, since $\comp{R}\otimes_R X$ and
$\comp{R}\otimes_R Y$ are isomorphic.
\end{proof}

Here is an example discovered by S.~Paul Smith to illustrate that, in the set-up of the
theorem above, $\gdim_{\vf}(X)$ and $\gdim_{\psi}(Y)$ need not be equal;  
see however~\cite[(8.2.4)]{avramov:homolhattfe}.

\begin{ex} \label{ex:not-equal}
  Let $R$ be a field, $S$ the localized polynomial ring $R[X]_{(X)}$, and let $T$ be a
  field extension of $R$, with  $\rank_RT=\rank_RS$; in particular, $S$ and $T$ are
  isomorphic as $R$-modules. Let $\vf\colon R\to S$ and $\psi\colon R\to T$ be the
  canonical inclusions.  Because $R$ is a field, both $\gdim_{\vf}(S)$ and
  $\gdim_{\psi}(T)$ are finite; see Theorem~\ref{thm:gor-char-1}.  By the
  Auslander-Bridger formula~\ref{thm:AB}, one has
\[ 
\gdim_{\vf}(S) = - 1\qquad\text{and}\qquad \gdim_{\psi}(T)= 0\,.
 \]
\end{ex}

The corollary below extends~\ref{subprop:4}. It applies, for instance, when $X$ is
homologically finite over $S$ and $\vf$ is module-finite.

\begin{cor} \label{cor:vf-is-finite}
  Let $\vf\colon R\to S$ be a local homomorphism and $X$ a complex of $S$-modules.  If
  $\HH(X)$ is finite over $R$, then
\[ 
\gdim_{\vf}(X) =\gdim_R(X). 
\]
\end{cor}

\begin{proof}
  Theorem~\ref{thm:indep-of-vf} applied to the homomorphisms $\vf$ and $\mathrm{id}_R$
  says that $\gdim_{\vf}(X)$ and $\gdim_R(X)$ are simultaneously finite.  When they are
  finite, the Auslander-Bridger formula and Lemma~\ref{lem:depth-eq} yield the desired
  equality.
\end{proof}

\section{Comparison with Gorenstein flat dimension}
\label{sec:gfd}
Keeping in mind the conclusions of the preceding section, and Proposition
\ref{prop:likeAF4.3}, it is natural to ask how G-dimension over $\vf$ compares with other
extensions of G-dimension to the non-finite arena. It turns out that the finiteness of
$\gdim_\vf(X)$ is equivalent to the finiteness of $\gfd_R(X)$, the \emph{G-flat dimension}
of $X$ over $R$, at least when $R$ has a dualizing complex. A more precise statement is
contained in Theorem \ref{thm:gdvsgfd} below; it is analogous to
Property~\ref{prop:pdvsfd} dealing with projective dimensions. We begin by recalling the
relevant definitions.

\begin{para}
\label{para:gfd}
Let $R$ be a commutative Noetherian ring. An $R$-module $G$ is said to be \emph{G-flat} if
there exists an exact complex of flat $R$-modules
\[ 
F= \cdots\xrightarrow{\partial_{i+1}}F_i\xrightarrow{\partial_{i}}F_{i-1}
\xrightarrow{\partial_{i-1}}\cdots
\] 
with $\coker(\partial_1)=G$ and $E\otimes_RF$ exact for each injective $R$-module $E$.
Note that any flat module is G-flat. Thus, each homologically bounded complex of
$R$-modules $X$ admits a G-flat resolution, and one can introduce its G-flat dimension to
be the number
\[
\gfd_R(X)\colon=\inf\{\sup\{n\mid G_n\ne 0\}\mid \text{$G$ a G-flat resolution of $X$}\}
\]
The reader may consult ~\cite{christensen:gd} or the book of Enochs and
Jenda~\cite{enochs:gf} for details.
\end{para}
 
Now we state one of the main theorems of this section; it implies Theorem D from the
introduction because when $R$ is a quotient of a Gorenstein ring, it has a dualizing
complex. As noted in the introduction, Foxby has derived the inequalities below assuming
only that the formal fibres of $R$ are Gorenstein. Also, the simultaneous finiteness of
$\gdim_{\vf}(S)$ and $\gfd_R(S)$ is \cite[(5.2)]{christensen:new}.  

\begin{thm}
\label{thm:gdvsgfd}
Suppose $R$ has a dualizing complex. Let $\vf\colon R\to S$ be a local homomorphism, and
$X$ a homologically finite complex of $S$-modules. Then
\[ 
\gfd_R(X) - \embdim(\vf) \leq \gdim_{\vf}(X) \leq \gfd_R(X)\,.
 \]
 In particular, $\gdim_{\vf}(X)$ is finite if and only if $\gfd_R(X)$ is finite.
\end{thm}

Observe that doing away with the hypothesis that $R$ has a dualizing complex would provide us
with another proof of Theorem \ref{thm:indep-of-vf}.  One can obtain useful bounds even
when $R$ has no dualizing complex; this is explained in \ref{para:closing},

\begin{para}
\label{proof:outline}
The proof calls for considerable preparation and is given in \ref{proof:gfdvsrfd}.  Here
are the key steps in our argument:
  
\emph{Step $1$.} We verify that $\gdim_{\vf}(X)$ and $\gfd_R(X)$ are simultaneously
finite.
This is an immediate consequence of \cite[(4.3)]{christensen:new} and
Theorem~\ref{prop:likeAF4.3}.
  
\emph{Step $2$.} We prove, in Theorem \ref{prop:gfd=rfd}, that if $\gfd_R(X)$ is finite,
then it coincides with the number $\rfd_R(X)$, whose definition is
recalled below.
This step constitutes the bulk of work in this section and builds on recent work of
Christensen, Frankild, and Holm~\cite{christensen:new}. They have informed us that they
can prove the same result by using the methods in \cite{holm:ghd}.

\emph{Step $3$.} The last step consists of verifying that when $\gdim_{\vf}(X)$ is finite,
it is sandwiched between $\rfd_R(X)-\embdim(\vf)$ and $\rfd_R(X)$.
The details of this step were worked out in conversations with Foxby, and we thank him for
permitting us to present them here.
\end{para}

\begin{para}
\label{para:rfd}
In the next few paragraphs, $R$ denotes a commutative Noetherian ring, not necessarily
local, and $W$ a homologically bounded complex of $R$-modules; we do not assume that
$\HH(W)$ is finite. The \emph{large restricted flat dimension} of $W$ over $R$, as
introduced in \cite{rfd}, is the quantity
\[
\rfd_R(W) = \sup\{\sup(F\lotimes_RW)\mid \text{$F$ an $R$-module with $\fdim_R(F)$
  finite}\}
\]
This number is finite, as long as $\HH(W)$ is nonzero and the Krull dimension of $R$ is
finite; see~\cite[(2.2)]{rfd}. It is useful to keep in mind an alternative formula
\cite[(2.4)]{rfd} for computing this invariant:
\begin{equation}
\rfd_R(W)=\sup\{\depth R_\p - \depth_{R_\p}(W_\p) \mid \p\in\spec R \}\,.
\end{equation}
\end{para}

We collect a few simple observations concerning this invariant.

\begin{lem}
\label{lem:rfd}
Let $\psi\colon R\to T$ and $\kappa\colon T\to T'$ be homomorphisms of commutative
Noetherian rings, and let $W$ and $Y$ be homologically bounded complexes of $R$-modules
and of $T$-modules respectively.
\begin{enumerate}[{\quad\rm(1)}]
\item If $\psi$ is faithfully flat, then
\[
\rfd_T(T\otimes_RW) = \rfd_R(W)\quad\text{and}\quad \rfd_T(Y) \geq \rfd_R(Y)\,.
\]
\item If $\kappa$ is faithfully flat, then \( \rfd_R(Y) = \rfd_R(T'\otimes_TY)\,.  \)
\end{enumerate}
 \end{lem}
\begin{proof}
  Let $F$ be an $R$-module and let $G$ be an $T$-module.
  
  Proof of (1). The flatness of $\psi$ implies
\begin{enumerate}[{\quad\rm(a)}]
\item if $\fdim_R(F)$ is finite, then so is $\fdim_T(F\otimes_RT)$;
\item if $\fdim_T(G)$ is finite, then so is $\fdim_R(G)$.
\end{enumerate}
Remark (a), combined with the isomorphisms
\[
(F\otimes_RT)\lotimes_T(T\otimes_RW) \simeq (F\otimes_RT)\lotimes_RW \simeq T\otimes_R
(F\lotimes_RW)
\]
and the faithful flatness of $\psi$, implies $\rfd_T(T\otimes_RW)\geq \rfd_R(W)$. The
opposite inequality follows from (b) and the associativity isomorphism
\[
G\lotimes_T(T\otimes_RW)\simeq G\lotimes_RW\,.
\]
This justifies the equality. The inequality is a consequence of (a) and the isomorphism
$(F\otimes_RT)\lotimes_TY \simeq F\lotimes_RY$.

As to (2): it is an immediate consequence of the isomorphism
\[
F\lotimes_R(T'\otimes_TY)\simeq (F\lotimes_RY)\otimes_TT'
\]
and the faithful flatness of $\kappa$.
\end{proof}

The next lemma gives a lower bound for the large restricted flat dimension.

\begin{lem}
\label{lem:rfd1}
If $\psi\colon R\to T$ is a local homomorphism and $Y$ is a complex $T$-modules, then
\[
\rfd_R(Y) \geq \depth R - \depth_T(Y)\,;
\]
equality holds if $Y$ is homologically finite over $R$ and $\gdim_R(Y)$ is finite.
 \end{lem}
\begin{proof}
  The inequality is a consequence of \ref{para:rfd}.1 and the (in)equalities
\[
\depth_R(Y) = \depth_T(\m T, Y) \leq \depth_T(Y)
\]
where the first one is by~\cite[(5.2.1)]{iyengar:dfcait} and the second is
by~\cite[(5.2.2)]{iyengar:dfcait}. If $Y$ is homologically finite over $R$ and
$\gdim_R(Y)$ is finite, then
\begin{align*}
  \depth R - \depth_TY&= \depth R - \depth_R Y \\
  &= \gdim_R(Y) \\
  &\geq \gdim_{R_\p}(Y_\p) \\
  &=\depth R_\p - \depth_{R_\p}(Y_\p)
 \end{align*}
 where the first equality is by Lemma \ref{lem:depth-eq}, the second and fourth are by the
 classical Auslander-Bridger formula, while the inequality is well known; see
 \cite[(2.3.11)]{christensen:gd}. In view of \ref{para:rfd}.1, this justifies the claimed
 equality.
\end{proof}

The next step towards Theorem \ref{thm:gdvsgfd} is the formula below. It may be viewed as
an Auslander-Buchsbaum formula for complexes of finite G-flat dimension, for is strikingly
similar to one for complexes of finite \emph{flat} dimension: $\depth_R(W) = \depth R -
\sup(k\lotimes_RW)$ when $\fdim_R(W)$ is finite; see \cite[(2.4)]{foxby:daafuc}.  What is
more, $E\lotimes_RW\simeq E\otimes_RG$, where is $G$ any G-flat resolution of $W$; this is
contained in \cite[(3.15)]{christensen:new}.

\begin{thm}
\label{prop:gfd-depth}
Let $(R,\m,k)$ be a local ring and $E$ the injective hull of $k$.  If $W$ is a complex of
$R$-modules with $\gfd_R(W)$ finite, then
\[
\depth_R(W) = \depth R - \sup(E\lotimes_RW)
\]
In particular, $\sup(E\lotimes_RW)$ is finite if and only if $\depth_R(W)$ is finite.
\end{thm}

\begin{proof}
  Let $\comp R$ denote the $\m$-adic completion of $R$. Faithful flatness of the
  completion homomorphism $R\to \comp R$ implies that each injective $\comp R$-module is
  injective also as an $R$-module. This remark and an elementary argument based on the
  definition of G-flat dimension entail: $\gfd_{\comp R}(\comp R\otimes_R W)\leq
  \gfd_R(W)$; see also Holm \cite[(3.10)]{holm:ghd}. Moreover
\[
\depth_{\comp R}(\comp R\otimes_RW) = \depth_RW \quad\text{and}\quad \depth {\comp
  R}=\depth R
\]
Finally, $E\lotimes_R W\simeq E\lotimes_{\comp R}(\comp R\otimes_RW)$, since $E$ has the
structure of an $\comp R$-module. Also, $E$ is the injective hull of $k$ over $\comp R$.
The upshot of this discussion is that one can replace $R$ and $W$ by $\comp R$ and $\comp
R\otimes_RW$, respectively, and assume that $R$ is complete.  In particular, $R$ has a
dualizing complex $D$.

The G-flat dimension of $W$ is finite, so it follows from \cite[(4.3)]{christensen:new}
that $W$ belongs to $\A(R)$, the Auslander category of $R$; see \ref{para:auslander}.
Thus, the canonical morphism $W\to \rhom_R(D,D\lotimes_RW)$ is an isomorphism, and this
starts the chain of isomorphisms
\begin{align*}
  \rhom_R(k,W) &\simeq \rhom_R(k,\rhom_R(D,D\lotimes_RW)) \\
  &\simeq \rhom_R(D\lotimes_Rk,D\lotimes_RW) \\
  &\simeq \rhom_k(D\lotimes_Rk,\rhom_R(k,D\lotimes_RW))\,.
\end{align*}
The second isomorphism is adjunction, so is the last one, since $D\lotimes_Rk$ is
isomorphic to a complex of vector spaces over $k$; this latter fact is clear once we
compute it with a free resolution of $D$.  The complex $D\lotimes_RW$ is homologically
bounded, since $W$ is in $\A(R)$, so the isomorphisms above
with~\cite[(1.5)]{foxby:daafuc} yield
\[
\sup(\rhom_R(k,W)) = \sup(\rhom_R(k,D\lotimes_RW)) - \inf(D\lotimes_Rk)\,.
\]
For each complex $X$ of $R$-modules, $\sup(\rhom_R(k,X))=-\depth_R(X)$ by
\cite[(2.1)]{foxby:daafuc}, and $\inf(D\lotimes_Rk)=\inf(D)$ since $D$ is homologically
finite, so the displayed equality translates to
\[
\depth_R(W) = \depth_R(D\lotimes_RW) + \inf(D)\,.
\]
Here is a crucial swindle: since $\gfd_R(W)$ is finite, so is $\gfd_R(\lch_\m(W))$, where
$\lch_\m(W)$ is the derived local cohomology of $W$ with respect to $\m$; this is by
\cite[(5.9)]{christensen:new}. Thus, the formula above applies to $\lch_\m(W)$ as well,
and reads
\[
\depth_R(\lch_\m(W)) = \depth_R(D\lotimes_R\lch_\m(W)) + \inf(D)\,.
\]
The homology modules of $\lch_\m(W)$ are all $\m$-torsion, so \cite[(2.7)]{foxby:daafuc}
yields the first the equality below, while \cite[(2.1)]{foxby:daafuc} provides the second
one
\[
\depth_R(\lch_\m(W))=-\sup(\lch_\m(W))=\depth_R(W)\,.
\]
Now, $\lch_\m(D) \simeq \shift^d E$ with $d=\inf(\lch_\m(D))$, where $\shift^d(-)$ denotes
a shift of $d$ steps to the left, so
\[
D\lotimes_R\lch_\m(W)\simeq \lch_\m(D)\lotimes_R W \simeq \shift^d (E\lotimes_R W)\,.
\]
The first isomorphism may be justified by invoking \cite[(3.1.2)]{lipman}.  The injective
hull $E$ is $\m$-torsion, so the homology modules of $E\lotimes_RW$ are $\m$-torsion:
compute via a free resolution of $W$. Thus, $\depth_R(E\lotimes_RW) = -\sup(E\lotimes_RW)$
by \cite[(2.7)]{foxby:daafuc}. Combining the preceding equalities gets us
\[
\depth_R(W) = -\sup(E\lotimes_RW) + \inf(\lch_\m(D)) + \inf(D)\,.
\]
Since $R$ itself has finite G-flat dimension, this formula with $R$ substituted for $W$
reads: \( \depth R = \inf(\lch_\m(D)) + \inf(D)\,.  \) Feeding this back into the formula
above completes the proof.
\end{proof}

In the case when $\HH(W)$ is concentrated in a single degree, the next theorem is part of
\cite[(3.19)]{holm:ghd}. 

\begin{thm}
\label{prop:gfd=rfd}
Let $R$ be a commutative Noetherian ring, and $W$ a complex of $R$-modules.  If
$\gfd_R(W)$ is finite, then $\gfd_R(W)=\rfd_R(W)$.
 \end{thm}

 \begin{proof}
   The proof is the following sequence of equalities:
\begin{align*}
  \gfd_R(W) &= \sup\{\sup(I\lotimes_RW)\mid \text{$I$ an injective $R$-module}\} \\
  &= \sup\{\sup(E(R/\p)\lotimes_R W) \mid \p\in\spec R\} \\
  &= \sup\{\sup(E(R/\p)\lotimes_{R_\p} W_\p) \mid \p\in\spec R\} \\
  &= \sup\{\depth R_\p - \depth_{R_\p}W_\p \mid \p\in\spec R\} \\
  &=\rfd_R(W)
\end{align*}
The first one is \cite[(2.8)]{christensen:new}; the second follows from this, given the
structure of injective modules over commutative Noetherian rings; the third equality is
due to the isomorphism $E(R/\p)\lotimes_R W \simeq E(R/\p)\lotimes_{R_\p} W_\p$, the
fourth is by Theorem \ref{prop:gfd-depth}, as $\gfd_{R_\p}(W_\p)\leq \gfd_R(W)$; see
\cite[(5.2.7)]{christensen:gd}, whilst the last equality is \ref{para:rfd}.1.
 \end{proof}
 
 We pause to record a corollary.

\begin{cor}
\label{cor:rfd-basechange}
Let $\psi\colon R\to T$ be a faithfully flat homomorphism of commutative, Noetherian
rings, and let $W$ be a complex of $R$-modules. If $\gfd_R(W)$ is finite, then
\[
\gfd_T(T\otimes_RW)=\gfd_R(T\otimes_RW) = \gfd_R(W)\,.
 \]
 \end{cor}
\begin{proof}
  When $\gfd_R(W)$ is finite, so are $\gfd_R(T\otimes_RW)$ and $\gfd_T(T\otimes_RW)$; the
  second by \cite[(3.10)]{holm:ghd} and the first follows from the easily verifiable
  remark: if $G$ is a G-flat $R$-module, the so is $F\otimes_RG$ for any flat $R$-module
  $F$. The desired equalities are now a consequence of Theorem \ref{prop:gfd=rfd} and
  Lemma \ref{lem:rfd}.
\end{proof}

This result prompts us to raise the


\begin{Question}
\label{quest}
Does the conclusion of Corollary \ref{cor:rfd-basechange} remain true without assuming
\emph{a priori} that $\gfd_R(W)$ is finite?
 \end{Question}

 Here, at last, is the proof of Theorem \ref{thm:gdvsgfd}; before jumping into it, the
 reader may wish to glance at \ref{proof:outline}, which outlines the basic argument.

\begin{para}\textit{Proof of Theorem~\ref{thm:gdvsgfd}.}
\label{proof:gfdvsrfd}
To begin with
\[
\gdim_{\vf}(X) < \infty \iff X \in A(R) \iff \gfd_R(X)<\infty\,,
\]
where the first biimplication is by Proposition~\ref{prop:likeAF4.3}, while the second one
is contained in~\cite[(4.3)]{christensen:new}.  Thus, one may assume that both
$\gdim_{\vf}(X)$ and $\gfd_R(X)$ are finite.  In this case, thanks to theorems
\ref{thm:AB} and \ref{prop:gfd=rfd}, what we need to prove is that
\begin{equation*}
\rfd_R(X) - \embdim(\vf) \leq \depth R - \depth_S X \leq \rfd_R(X) 
\tag{$\dagger$}
\end{equation*}
when $\gdim_{\vf}(X)$ is finite.  The inequality on the right is contained in Lemma
\ref{lem:rfd1}. That leaves us with the one on the left.

Let $\comp S$ be the completion of $S$ at its maximal ideal and set $\comp X =\comp
S\otimes_S X$.  By Lemma \ref{lem:rfd}.2, the faithful flatness of the homomorphism $S\to
\comp S$ implies $\rfd_R(\comp X) = \rfd_R(X)$. The other quantities involved in
($\dagger$) also remain unchanged if we substitute $\comp S$ for $S$ and $\comp X$ for
$X$, so we may do so and thereby assume that $S$ is complete. With $R\to R'\to S$ a
minimal Cohen factorization of $\vf$, Lemma \ref{lem:rfd} provides the inequality below
\begin{align*}
  \rfd_R(X) &\leq \rfd_{R'}(X)  \\
  & = \depth R' - \depth_{R'}(X)  \\
  &=\depth R + \embdim(\vf) - \depth_{R'}(X) \\
  &=\depth R + \embdim(\vf) - \depth_S(X)
\end{align*}
Lemma \ref{lem:rfd1} explains the first equality; the second holds as $R\to R'$ is flat
and $R'/\m R'$ is regular, and the last holds because $R'\to S$ is surjective. \qed
\end{para}

\begin{para}
\label{para:closing}
Let $\vf\colon (R,\m,k) \to S$ be a local homomorphism and $X$ a homologically finite
complex of $S$-modules. Let $\comp R$ denote the $\m$-adic completion of $R$, and $\wt S$
the $\m S$-adic completion of $S$. Since $\comp R$ has dualizing complex, it follows from
\ref{subprop:1} and Theorem \ref{thm:gdvsgfd} that
\[
\gfd_{\comp R}(\wt S\otimes_SX) - \embdim(\vf) \leq \gdim_{\vf}(X) \leq \gfd_{\comp R}(\wt
S\otimes_SX)\,.
\]
 
At any rate, one has the consolation of knowing a partial result:

\begin{prop}
  Let $\vf\colon R\to S$ be a local homomorphism. Each homologically finite complex of
  $S$-modules $X$ satisfies $\gdim_{\vf}(X)\leq \gfd_R(X)$.
\end{prop}
\begin{proof}
  The plan is to reduce to the case where $R$ is complete and then apply Theorem
  \ref{thm:gdvsgfd}; confer the proof of Theorem \ref{thm:indep-of-vf}. Let $K$ be the
  Koszul complex on minimal set of generators for $\m$, the maximal ideal of $R$. Thus,
  $\pdim_R(K)=\embdim R = \pdim_S(K\otimes_RS)$.  Now, if $G$ is a G-flat resolution of
  $X$ over $R$, then $K\otimes_RG$ is a G-flat resolution of $K\otimes_RX$. This implies
  that
\[
\gfd_R(K\otimes_RX) \leq \gfd_R(X) + \embdim(R)\,.
\]
Moreover, since $K\otimes_RX\cong (K\otimes_RS)\otimes_SX$, Theorem
\ref{thm:stability-tensor} applied to the diagram $R\to S\xrightarrow{=}S$, and with
$P=(K\otimes_RS)$, yields
\[
\gdim_{\vf}(K\otimes_RX) = \gdim_{\vf}(X) + \embdim(R)\,.
\]
Thus, it suffices to prove the desired inequality for the complex of $S$-modules
$K\otimes_RX$; in particular, one may pass to $K\otimes_RX$ and assume $\m\cdot \HH(X)=0$.
Now, we adopt the notation of \ref{para:closing}, where we noted that
\[
\gdim_{\vf}(X) \leq \gfd_{\comp R}(\wt S\otimes_SX)\,.
\]
It is elementary to verify that the canonical homomorphism of complexes of $\comp
R$-modules $\comp R\otimes_R X \to \wt S\otimes_S X$ is a homology isomorphism, since
$\m\cdot\HH(X)=0$. This gives us the equality below:
\[
\gfd_{\comp R}(\wt S\otimes_SX)=\gfd_{\comp R}(\comp R\otimes_RX)\leq \gfd_R(X)\,;
\]
the inequality is the version for complexes of \cite[(3.10)]{holm:ghd}, and may be deduced
directly from the definitions. To complete the proof, put together the composed inequality
above with the penultimate one.
\end{proof}

This proposition leads to analogues of the theorems in Section \ref{sec:apps}, with
$\gfd_R(-)$ playing the role of $\gdim_{\vf}(-)$.  The result below, which parallels
Theorem \ref{thm:gdim(f)-finite}, is one such; in it, for any complex of $R$-modules, we
write ${}^{\vf^n}\!X$ to indicate that $R$ acts on $X$ via $\vf^n$.

\begin{thm} 
  Let $\vf\colon (R,\m)\to (R,\m)$ be a local homomorphism such that
  $\vf^i(\m)\subseteq\m^2$ for some integer $i\geq 1$.  The following conditions are
  equivalent.
\begin{enumerate}[{\quad\rm(a)}] 
\item The ring $R$ is Gorenstein.
\item $\gfd_R({}^{\vf^n}\!R)$ is finite for each integer $n\geq 1$.
\item There is a homologically finite complex $P$ of $R$-modules with $\pdim_R(P)$ finite
  and $\gfd_R({}^{\vf^n}\!P)$ finite, for some integer $n\geq 1$.
\end{enumerate}
\end{thm}
\begin{proof}
  Over a Gorenstein ring, any module has finite G-flat dimension; see
  \cite[(5.2.10)]{christensen:gd}.  This justifies ``(a) $\implies$ (b)'', while ``(b)
  $\implies$ (c)'' is trivial.
  
  ``(c) $\implies$ (a)''.  The preceding proposition yields that $\gdim_{\vf^n}(P)$ is
  finite, so it remains to invoke the corresponding implication in Theorem
  \ref{thm:gdim(f)-finite}.
\end{proof}
\end{para}

\begin{para}
  The hypotheses of the preceding result are satisfied when $\vf$ is the Frobenius
  endomorphism of a local ring of $R$ of positive prime characteristic. In this case, one
  can add a fourth equivalent condition to those given above:
\begin{enumerate}[{\quad\rm(a)}]
\item[{\quad\rm(b$'$)}] the $R$-module ${}^{\vf^n}\!R$ is G-flat for one integer $n\geq
  1$.
\end{enumerate}
Indeed, it is clear from \ref{para:rfd}.1 that $\rfd_R({}^{\vf^n}\!R)=0$ for each integer
$n\geq0$. Therefore, by \cite[(3.19)]{holm:ghd}, or by its successor, Theorem
\ref{prop:gfd=rfd}, one obtains that the $R$-module ${}^{\vf^n}\!R$ has finite G-flat
dimension if and only if it is already G-flat.
\end{para}

\section*{Acknowledgments}
While this article was being written, the first author was collaborating with Luchezar
Avramov and Claudia Miller on \cite{avramov:homolhattfe}, another study of local
homomorphisms.  We are grateful to them for numerous discussions.
  
We thank Lars Christensen, Anders Frankild, and Henrik Holm for making
\cite{christensen:new} available to us and for correspondence concerning Section 8.  Our
work in that section owes an intellectual debt to Hans-Bj\o rn Foxby; we thank him for
comments and criticisms, and for generously sharing his ideas.

Thanks are also due to S.~Paul Smith, Diana White, Yongwei Yau, and the referee for useful
suggestions.

\end{document}